\documentclass[a4paper]{article}
\usepackage[utf8]{inputenc}
\usepackage[T1]{fontenc}
\usepackage{tikz}
\usepackage{graphicx}
\usepackage{bm}
\usepackage{amsmath}
\usepackage{xcolor}
\usepackage[backend=bibtex]{biblatex}%

\usepackage{authblk}
\usepackage[version=4]{mhchem}
\usepackage[section]{placeins}
\usepackage{pgfplots}
\usepackage{subcaption}
\usepackage{xcolor}
\usepackage{placeins}
\usepackage{hyperref}
\usepackage[normalem]{ulem}

\definecolor{pltc0}{HTML}{1f77b4}
\definecolor{pltc1}{HTML}{ff7f0e}
\definecolor{pltc2}{HTML}{2ca02c}
\definecolor{pltc3}{HTML}{d62728}
\definecolor{pltc4}{HTML}{9467bd}
\definecolor{pltc5}{HTML}{8c564b}
\definecolor{pltc6}{HTML}{e377c2}
\definecolor{pltc6}{HTML}{7f7f7f}
\definecolor{pltc8}{HTML}{bcbd22}
\definecolor{pltc9}{HTML}{17becf}

\usepackage{comment}
\usepackage{amsmath}
\usepackage{amssymb}       
\usepackage{amsthm}      

\usepackage{color}
\usepackage{bm,bbm}
\usepackage{wrapfig}
\usepackage{overpic}
\usepackage{siunitx}
\usepackage{soul}
\usepackage{enumitem}

\theoremstyle{plain}

\theoremstyle{definition}

\theoremstyle{plain}

\definecolor{MyDarkGreen}{rgb}{0,0.45,0}

\def\trait #1 #2 #3 {\vrule width #1pt height #2pt depth #3pt}
\def\fin{\hfill
        \trait .3 5 0
        \trait 5 .3 0
        \kern-5pt
        \trait 5 5 -4.7
        \trait 0.3 5 0
\medskip}

\newif\ifANKI

\newcounter{numbs}
\newcounter{numbi}
\newcounter{numbii}

\ANKIfalse

\ifANKI

\else

\fi

\newcommand{\PGRAPH}[1]{\noindent\textbf{#1}}

\newcommand{\REAL}{\mathbbm{R}}

\newcommand{\av}{\mathbf{a}}
\newcommand{\bv}{\mathbf{b}}
\newcommand{\ccv}{\mathbf{c}}

\newcommand{\rv}{\mathbf{r}}

\newcommand{\xv}{\mathbf{x}}

\newcommand{\as}{a}
\newcommand{\bs}{b}
\newcommand{\cs}{c}

\newcommand{\fs}{f}
\newcommand{\gs}{g}

\newcommand{\is}{i}
\newcommand{\js}{j}

\newcommand{\ns}{n}

\newcommand{\rs}{r}

\newcommand{\ts}{t}

\newcommand{\xs}{x}

\newcommand{\As}{A}
\newcommand{\Bs}{B}
\newcommand{\Cs}{C}

\newcommand{\Ls}{L}

\newcommand{\Ps}{P}

\newcommand{\Ts}{T}
\newcommand{\Us}{U}
\newcommand{\Vs}{V}

\newcommand{\fss}[1]{f_{#1}}

\newcommand{\iss}[1]{i_{#1}}
\newcommand{\jss}[1]{j_{#1}}

\newcommand{\nss}[1]{n_{#1}}

\newcommand{\rss}[1]{r_{#1}}

\newcommand{\xss}[1]{x_{#1}}

\newcommand{\Ass}[1]{A_{#1}}

\newcommand{\Gss}[1]{G_{#1}}

\newcommand{\Lss}[1]{L_{#1}}

\newcommand{\Pss}[1]{P_{#1}}

\newcommand{\Psz}[2]{\Ps_{#1}^{#2}}

\newcommand{\calF}{\mathcal{F}}

\newcommand{\calL}{\mathcal{L}}

\newcommand{\calO}{\mathcal{O}}

\newcommand{\tA}{\mathcal{A}}
\newcommand{\tB}{\mathcal{B}}
\newcommand{\tC}{\mathcal{C}}

\newcommand{\tF}{\mathcal{F}}

\newcommand{\tATT}{\mathcal{A}^{TT}}

\newcommand{\DIM} {d}              

\newcommand{\nlen}{\hspace{-0.2mm}}

\newcommand{\norm}   [2]{|\nlen|#1|\nlen|_{#2}}

\newcommand{\abs}    [1]{|#1|}

\newcommand{\EOD}{\end{document}}

\newcommand{\roundPrecision}{2}
\newcommand{\TT}{TT}

\newcommand{\betav}{{\bm\beta}}   

\newcommand{\fsn}{\fs_{\ns}}

\newcommand{\range}[1]{]#1]\phantom{[[} }

\newcommand{\RNDG}{\scalebox{0.9}{{\sf rndg}}}
\newcommand{\TTSVD}{\scalebox{0.9}{{\sf TT-SVD}}}
\newcommand{\cTT}{\mathbbm{TT}}
\newcommand{\toll}{\epsilon}

\newcommand{\tAss}[1]{\tA_{#1}}

\newcommand{\MAXVOL}{\scalebox{0.9}{{\sf MaxVol}}}

\usepackage[a4paper, top=2cm, bottom=2cm, left=2.5cm, right=2.5cm]{geometry}

\title{A Fast, Accurate and Oscillation-free Spectral Collocation Solver for High-dimensional Transport Problems}

\date{}

\author[1]{Nicola Cavallini}
\author[2]{Gianmarco Manzini}
\author[3,4]{Daniele Funaro}
\author[1,2,*]{Andrea Favalli}

\affil[1]{European Commission, Joint Research Centre, Via Enrico Fermi, Ispra, I -21027 (Va), Italy}
\affil[2]{Los Alamos National Laboratory, P.O. Box 1663, Los Alamos, NM 87545, USA}
\affil[3]{Dipartimento di Scienze Chimiche e Geologiche, Universit\`a di Modena e Reggio Emilia, Via Campi 103, \phantom{${}^{3}$}41125 Modena, Italy}
\affil[4]{Istituto di Matematica e Applicata e Tecnologie Informatiche del CNR, via Ferrata 1, 27100 Pavia, Italy}

\affil[*]{andrea.favalli@ec.europa.eu}

\newcommand{\mymethod}{T$^2$S$^2$}
\newcommand{\mymethodextended}{\textbf{Tensor Train Superconsistent Spectral} }
\newcommand{\advection}{$\varepsilon\slash{\abs{\betav\Ls}}\ll1$}
\newcommand{\AdvRatio}{$\varepsilon\slash{\abs{\betav\Ls}}$}

\begin{document}
\maketitle

\begin{abstract}


Transport phenomena---describing the movement of particles, energy,
or other physical quantities---are fundamental in various scientific disciplines, including nuclear physics, plasma
physics, astrophysics, engineering, and the natural sciences.
However, solving the associated seven-dimensional transport equations poses a significant
computational challenge due to the curse of dimensionality.
We introduce the Tensor Train Superconsistent Spectral (\mymethod{})
solver to address this challenge, integrating Spectral Collocation for
exponential convergence, Superconsistency for stabilization in
transport-dominated regimes, and Tensor Train format for substantial
data compression.
\mymethod{} enforces a dimension-wise superconsistent condition
compatible with tensor structures, achieving extremely low compression
ratios, such as $\mathcal{O}(10^{-12})$, while preserving spectral
accuracy.
Numerical experiments on linear problems demonstrate that \mymethod{}
can solve high-dimensional transport problems in minutes on standard
hardware, making previously intractable problems computationally
feasible.
This advancement opens new avenues for efficiently and
accurately modeling complex transport phenomena.



\end{abstract}


\newcommand{\marco}[2]{\textcolor{red}{\sout{#1}~}\textcolor{magenta}{#2}~}
\newcommand{\TIMEMARCHING}{second-order Crank-Nicolson~}

\newcommand{\PLAINSC}{\textrm{Plain-SC}}

\section*{Introduction}
Transport phenomena are fundamental across a wide 
range of scientific disciplines, including nuclear physics, 
plasma physics, high-energy density physics, astrophysics, 
condensed matter physics, atmospheric science, oceanography, and engineering. 
These processes describe the movement and interaction of particles, energy, 
and other physical quantities within various media.  
At the heart of modeling these phenomena are transport equations, 
which comprehensively describe how quantities propagate and interact. 
For example, the neutron transport equation is fundamental in applications 
such as reactor design and radiation detection instrumentation, 
while the Vlasov equation is pivotal in plasma physics, 
describing the evolution of charged particles under electromagnetic fields, 
and is fundamental in applications from fusion to space science. 
\cite{Zank:2014,Williams:2013,Jungel:2009,Cercignani:2012,Bell-Glasstone:1970,Graziani:2006}
Although the precise nature of transport phenomena 
may vary across applications, the underlying mathematical structure 
remains consistent. A significant challenge in solving transport equations is the 
associated computational cost, given that a general solution spans seven dimensions:  
three spatial coordinates, two angles, one energy (or speed), and one time. 
As the number of dimensions increases, the memory and computational resources 
required scale exponentially, limiting the feasibility of 
traditional numerical methods. The challenge has been aptly 
described by F. Graziani and G. Olson as ``Conquering the Seven-Dimensional Mountain''
\cite{Graziani:2003}.

To establish a robust computational framework, we focus 
on a fundamental yet representative model of transport phenomena: 
the linear, time-dependent, convection-diffusion-reaction 
Partial Differential Equation (PDE). We consider the PDE:
\begin{align}
  \frac{\partial\fs}{\partial\ts} -\varepsilon\Delta_{\mathbf{x}}\fs
  +\betav\cdot\nabla_{\mathbf{x}}\fs +\rho\fs = \bs,
  \qquad\forall(\xv,\ts)\in\Omega\times\range{0,T},
  \label{eq:pblm:time:A}
\end{align}
The equation describes the evolution of the distribution 
function $f(\mathbf{x},t)$ in the seven-dimensional phase space 
$\Omega\times[0,T]$, with 
$\mathbf{x} = (x_{1},x_{2},\ldots,x_{d})^T\in\Omega\subset\mathbb{R}^{d}$, $d = 6$. 
Time is represented by $t\times [0,T]$ and $x_1, x_2, x_3$ refer to  the 3-D 
spatial variables. Depending on the specific form of equation (\ref{eq:pblm:time:A})  
$x_4, x_5, x_6$ either describe the momentum vector or 
$x_4, x_5$ represent the 2-D angular direction of flight, 
and $x_6$ the energy or speed. The computational domain 
$\Omega\subset\mathbb{R}^{d}$ is the 
$d$-dimensional, open hypercube with characteristic edge length $L$ 
and five-dimensional boundary $\Gamma$.
To have a mathematically well-posed model, 
we assume that $f$ satisfies Dirichlet boundary 
conditions on the domain boundary 
\begin{align}
  f = g\quad\textrm{on}~\Gamma\times[0,T],
  \label{eq:pblm:time:B}    
\end{align}
and an initial condition at $t=0$
\begin{align}
  \fs(\cdot,0)= f_0\quad\textrm{in}~\Omega.
  \label{eq:pblm:time:C}
\end{align}
In Eq.~\eqref{eq:pblm:time:A}, $\varepsilon$ is the diffusion 
coefficient, $\bm{\beta}$ is the convection field, $\rho$ 
is the reaction field, and $b$ is the right-hand side forcing term.
As a matter of notation, in the latter we drop the $\xv$ subscript 
in the differential operators notation: $\nabla_{\xv}$ becomes 
$\nabla$ and $\Delta_{\xv}$ becomes $\Delta$.

The numerical solution of high-dimensional problems 
in the framework of Partial Differential Equations (PDEs) 
faces a fundamental challenge often referred to as
the \emph{curse of dimensionality}~\cite{Bellman:1957}, 
where both memory requirements and computational costs 
scale exponentially with the number of dimensions. 
This is a serious limit for the effectiveness of numerical methods. 
While recent approaches have attempted to address this challenge 
through GPU acceleration \cite{Einkemmer:Moriggl:2021} 
and machine learning methods \cite{Zhu:2023,Han:Jentzen:2020}, 
these solutions introduce  further complexities, as they require 
either architecture-specific programming or extensive training data.

To address these computational challenges, we take a different path 
based on a three-level approach. 
First, we operate in the context of the 
\emph{low-rank approximation}\cite{Einkemmer:2024}. 
We adopt the Tensor Train format \cite{Oseledets:2011}, 
a linear algebra framework capable of efficiently handling 
high-dimensional data and mathematical operators, 
organized on tensor-product grids. Considerable effort has been 
devoted over the years to solving PDEs in Tensor Train format. 
Dolgov and coauthors tackled the three-dimensional Fokker-Planck 
equation  \cite{Dolgov-Khoromskij-Oseledets:2012}. 
Basic numerical schemes, such as uniform finite differences applied to 
the high-dimensional Poisson equation, are commonly used as testbed for the development of the Tensor Train linear 
system solvers \cite{Dolgov-Savostyanov:2014,Dolgov-Savostyanov:2015}. 
A hybrid approach for parametric PDEs in two dimensions 
can be found in \cite{Dolgov-Scheichl:2019}. In the context of 
semi-Lagrangian solvers, where special care has to be devoted to 
conservation of mass, energy, or momentum, Kormann \cite{Kormann:2015} 
proposed a Tensor Train solver for Vlasov Equations. 

Second, we numerically discretize 
the problem (\ref{eq:pblm:time:A}) using a variant of Spectral Collocation methods. 
The term ``collocation’’ indicates that the continuous model, i.e. 
the PDE, is evaluated, ``collocated’’, 
at specific grid points. These techniques are known for their long-recognized exceptional \textit{accuracy}
\cite{Canuto-Hussaini-Quarteroni-Zang:2007,Hesthaven-Gottlieb-Gottlieb:2007},  
commonly exhibiting an exponential rate of convergence for smooth solutions. 
This rate of convergence is also addressed as spectral accuracy. 
However, their application to higher-dimensional problems, 
in standard algebra format, has been constrained by the dense structure of the 
discrete operators and the complexity of solving the resulting linear 
systems \cite{Canuto-Hussaini-Quarteroni-Zang:2007,Funaro:1997}. 

Third, these methods face a fundamental trade-off between 
\textit{accuracy} and \textit{stability} \cite{Arnold:2018,Fish-Belytschko:2007}. 
This becomes particularly problematic in regimes dominated 
by first-order derivatives (\advection{}), where solutions present 
spurious oscillations that make the result unreliable. 
Funaro's concept of superconsistency \cite{Funaro:2002} 
cures such oscillations by performing collocation on modified nodes. 
The new position of the nodes is a function of the ratio 
\AdvRatio{} 
and the grid density \cite{Funaro:1997,DelIsle-Owens:2021,DelIsle-Owens:2023}. 
We should note that superconsistency alone does not reduce the computational costs; 
the solution of high-dimensional problems remains prohibitive.

\medskip
Our approach, 
\mymethodextended{} (\mymethod{}), 
builds on the three techniques described above. 
It addresses the limitations of each individual block through a 
novel dimension-wise application of the superconsistent method. 
This approach fundamentally redefines the original superconsistent framework, 
making it compatible with the Tensor Train format and enabling 
an efficient storage of the mathematical operators.
The method preserves the hallmark exponential accuracy of 
Spectral Collocation methods, while maintaining the stabilization 
properties with respect to the \AdvRatio{} ratio. 
At the same time \mymethod{} achieves extremely low compression 
ratios, such as $\mathcal{O}(10^{-12})$, as a result, it can 
solve high-dimensional transport problems in minutes on standard hardware, 
rendering previously intractable problems computationally feasible.
\medskip

The paper is organized as follows. 
The first section, ``Overview of the \mymethod{} Approach'', 
presents the method, followed by 
the second section, ``Results'', which demonstrates the validation 
of the method through numerical experiments. 
The third section, ``Discussion'', summarizes the method's strengths, 
applications and future directions. 
The fourth section, ``Methods'', provides the detailed mathematical 
framework underlying the \mymethod{} solver.

\section*{Overview of the \mymethod{} Approach}
\label{sec2:overview}
\begin{figure}[!t]
  \begin{center}
    \includegraphics[width=0.8\textwidth]{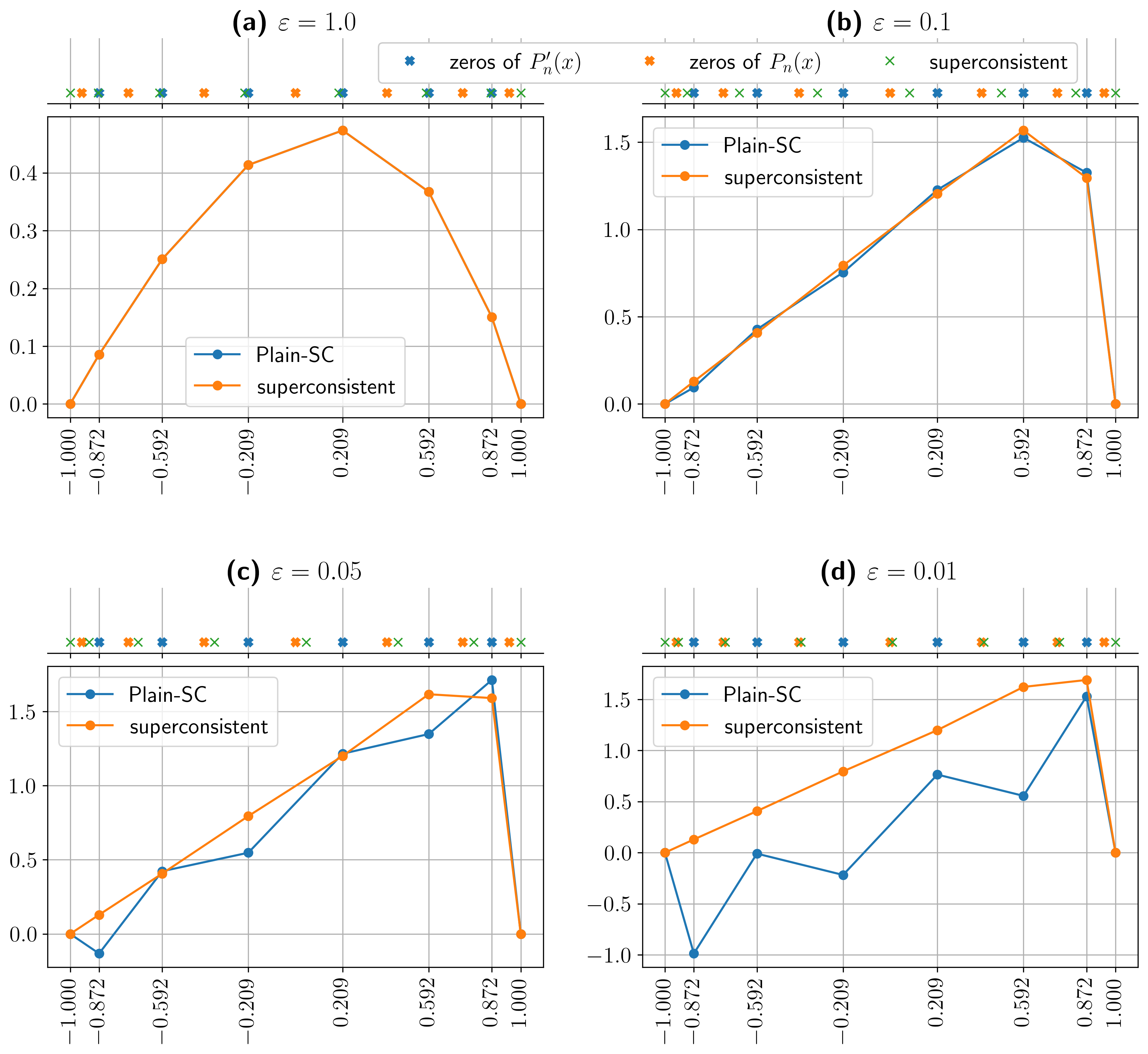}
  \end{center}
  \caption{Superconsistent approximation of Eq.\eqref{eq:adv-diff-1d}
    for the polynomial degree $n=7$, and $\beta=1$.
    On the top horizontal axis, the blue crosses mark the zeros of
    $P'_n$, optimal for diffusion; the orange crosses mark the zeros
    of $P_n$, optimal for convection; the green crosses mark the
    superconsistent collocation points for mixed regimes.
    Panels \textbf{(a)} and \textbf{(b)} show relatively accurate
    solutions when convection and diffusion are balanced, in particular
    in panel \textbf{(a)} the two solutions almost coincide.
    As convection dominates (panels \textbf{(c)} and \textbf{(d)}), the
    \PLAINSC{} solution is unreliable, whereas the
    superconsistent one does not present problems.  }
  \label{fig:superc-1d}
\end{figure}

This section outlines the key components of our methodology and
demonstrates how they synergize to create our computational framework.
We illustrate the underlying principles through one-
and two-dimensional examples,
while the detailed mathematical framework is provided in
Section ``Methods''.

\medskip
\PGRAPH{$(i)$ \emph{Spectral collocation methods.}}
Spectral collocation relies on two fundamental elements: \emph{basis
functions} and \emph{representation points}.
We employ Lagrange basis functions defined on $n+1$ representation
points \cite{Funaro:1997}, which are the $n-1$ zeros of
$\Pss{n}^{\prime}$, the derivative of the $n$-th degree Legendre
polynomial $P_n$, complemented with the two extrema of the interval
$[-1, 1]$.
We will refer this approach as \emph{plain spectral collocation}
(\PLAINSC{}) method.
The solution to a one-dimensional transport problem is a linear
combination of these basis functions.
The coefficients are determined by solving the linear system obtained
by \emph{collocating the equation} at the same representation points,
i.e., by evaluating the left-hand side and the right-hand side of
Equation~\eqref{eq:pblm:time:A} at such locations.
For nomenclature clarity, in one dimension, the number of degrees of
freedom ($\#\mathsf{dofs}$) equals $n+1$.
The multidimensional extension of the method is achieved through the
tensor product of the one-dimensional discretizations along all the
problem's directions.

\medskip
\PGRAPH{$(ii)$ \emph{Superconsistent formulation.}}  To illustrate the
concept of superconsistency we consider the time-independent
one-dimensional model

\begin{align}
  \beta \fs'(\xs) - \varepsilon \fs''(\xs) = 1 \quad x\ \mathrm{in~}\big[-1,1\big], \quad \fs(-1) = \fs(1) = 0,
  \label{eq:adv-diff-1d}
\end{align}

where $\beta>0$ and $\varepsilon>0$ are the scalar convective and
diffusion coefficients, respectively.
The \PLAINSC{} method, employing the zeros of $\Psz{n}{\prime}$, (blue
dots on the top axis of Figure~\ref{fig:superc-1d}), develops spurious
oscillations when the convective term becomes dominant relatively to
the diffusive one, meaning that \advection{}.
This effect is clearly visible in panels \textbf{(c)} and \textbf{(d)}
of Figure~\ref{fig:superc-1d}.

Funaro~\cite{Funaro:2002} cured this phenomenon by introducing the
superconsistent condition.
The superconsistent condition takes into account the differential
operator and is satisfied by introducing a different set of
\emph{collocation points} where the partial differential equation and
the right-hand side are evaluated.

According to the superconsistent condition, the zeros of $\Pss{n}'$
are optimal for the approximation of the diffusion operator
$\varepsilon f''$, while the zeros of $\Pss{n}$ are optimal for $\beta
f'$.
In hybrid regimes, the superconsistent points establish an optimal
intermediate distribution between these two configurations. As
illustrated in Figure~\ref{fig:superc-1d}, when convection and
diffusion terms are balanced (Figure~\ref{fig:superc-1d}\textbf{(a)}),
the superconsistent points are close to the zeros of $P'_n$.
As the diffusion coefficient $\varepsilon$ decreases maintaining $n$
fixed (Figure~\ref{fig:superc-1d}\textbf{(b)} and \textbf{(c)}), the
superconsistent points shift towards the zeros of $\Pss{n}$.
In the convection-dominated regime
(Figure~\ref{fig:superc-1d}\textbf{(d)}), the superconsistent points
nearly coincide with the zeros of $\Pss{n}$.

In summary, the superconsistent approach is 
based on two grids: the \emph{representation grid} and the 
actual \emph{collocation grid}.
The position of the collocation grid depends on the ratio \AdvRatio{}
and $n$.
Enforcing the superconsistent condition provides the numerical
strategy to locate the collocation points in order to achieve
stabilization with respect to \AdvRatio{}.
We use the representation points to build the Lagrange basis
functions.
These basis functions are used to evaluate the discrete operators at
the collocation grid.
The result of this \emph{collocation} operation, together with
boundary conditions, leads to the construction of the final linear
system.

\FloatBarrier

\medskip
\PGRAPH{$(iii)$ \emph{Tensor-train format.}}
Traditional spectral methods struggle with the \emph{curse of
dimensionality}~\cite{Bellman:1957}, where computational complexity
increases by scaling exponentially with the number of dimensions.
\begin{figure}[!t]
  \centering
  \includegraphics[width=1.\textwidth]{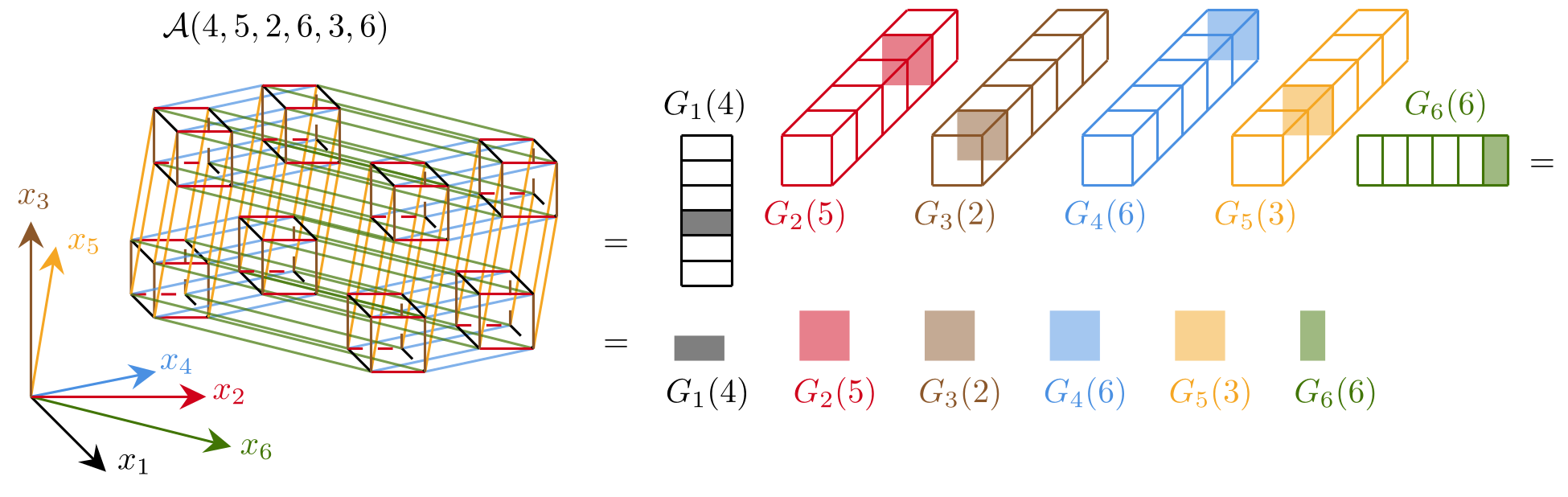}
  \caption{Focusing on the six spatial dimensions, a six-dimensional
    array $\mathcal{A}$ has size $\ns^6$, with $x_1, x_2, \ldots, x_6$
    coordinate axis of the six-dimensional space.
    The TT format represents the whole array on the left as the
    product of six cores on the right.
    The $\ell$-th core for $\ell=1,2,\ldots,6$ has size
    $\rss{\ell-1}\times\ns\times\rss{\ell}$, with $\rss{0}=\rss{6}=1$,
    so the first and last cores are matrices.
    ~The $(4,5,2,6,3,8)$-th entry of the full tensor is represented in
    the TT-format as the product
    $\mathcal{A}(4,5,2,6,3,8)=\Gss{1}(4)\Gss{2}(5)\Gss{3}(2)\Gss{4}(6)\Gss{5}(3)\Gss{6}(8)$,
    i.e., as the product of the $4$-th, $5$-th, $2$-nd, $6$-th,
    $3$-rd, $8$-th slices of the corresponding cores.  }
  \label{fig:TT-format}
\end{figure}
The Tensor-train (TT) format\cite{Oseledets:2011} efficiently represents a
$\DIM$-dimensional tensor as the product of $\DIM$ three-dimensional
arrays, the \emph{TT cores}, significantly reducing the memory
requirement.
Each tensor-train core, denoted by $\Gss{\ell}$ and indexed by
$\ell=1,2\ldots,\DIM$, has dimensions
$\rss{\ell-1}\times\nss{\ell}\times\rss{\ell}$, where $\nss{\ell}$ is
the original dimension size also called the \emph{mode size}, and
$\rss{\ell}$ are the \emph{TT-ranks} (with $\rss{0}=\rss{\DIM}=1$).
As shown in Figure~\ref{fig:TT-format} for a six-dimensional tensor,
an element of the multidimensional array
$\mathcal{A}\in\REAL^{\nss{1}\times\ldots\nss{6}}$ is reconstructed by
multiplying the corresponding core slices
$\Gss{\ell}\in\REAL^{\rss{\ell-1}\times\nss{\ell}\times\rss{\ell}}$,
$\ell=1,2,\ldots,6$.
The TT decomposition is not unique and is most effective when the
TT-ranks are much smaller than the mode size, i.e., for
$\rss{\ell}\ll\nss{\ell}$.
In such a case, this method efficiently compresses large
multidimensional arrays and, when applied to operators, mitigates the
high bandwidth issue arising from tensor products.

\medskip
We build the \mymethod{} solver combining 
principles $(i)$, $(ii)$, and $(iii)$
to extend superconsistency to high-dimensional
equations.
Our approach implements the superconsistent condition independently
along each dimension, treating them as decoupled problems.
This formulation works with two distinct time integration frameworks:
the \emph{unified space-time
formulation}\cite{tal1986spectral} that
handles time as an additional dimension, and the \emph{method-of-lines
formulation}\cite{Schiesser:1991,Hundsdorfer-Verwer:2003} that uses
conventional time-stepping schemes.
For the latter, we implement the Crank-Nicholson and the backward
Euler schemes.
Highly accurate variants of these schemes, such as
\cite{Ruuth-Spiteri:2003} are combinations of the mentioned methods.
Superconsistency is applied only to the six spatial dimensions,
demonstrating an oscillatory-free performance with respect to the
ratio \AdvRatio{}.

\begin{figure}[t!]
  \centering
  \includegraphics[width=1.\textwidth]{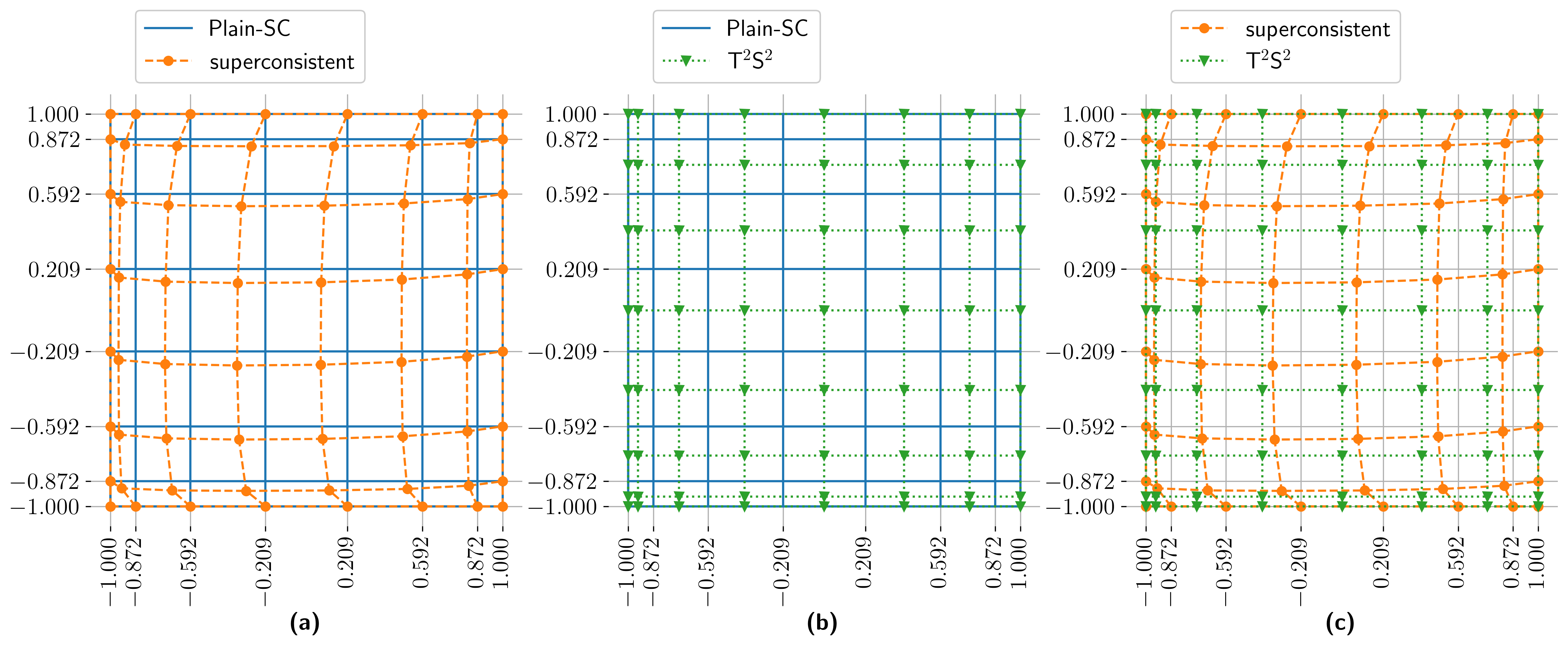}
  \caption{Left panel: \PLAINSC{} nodes and Funaro's superconsistent
    nodes\cite{Funaro:1997}.
    Such a grid has been evaluated for problem~\eqref{eq:adv-diff-1d}
    with parameters $\varepsilon={1}\slash{50}$, $\beta=(1,0.5)^T$ and
    $n=7$.
    Middle panel: \PLAINSC{} nodes and \mymethod{} collocation nodes,
    which have been independently computed in each direction.
    Right panel: Comparison between the 2D Funaro's superconsistent
    nodes\cite{Funaro:1997} and \mymethod{} collocation nodes.
  }
  \label{fig:superc-2d}
\end{figure}

In Figure~\ref{fig:superc-2d} we show a two dimensional grid 
obtained using \mymethod{}.
The resulting Cartesian structure naturally fits with the tensor-train
format and extends straightforwardly to higher dimensions.
This approach differs from Funaro's original
formulation\cite{Funaro:1997} , which applies superconsistency to the
full multidimensional problem and produces non-Cartesian grid
deformations, see Fig.~\ref{fig:superc-2d}.
Our method is stable across physically relevant parameter ranges and
preserves the exponential accuracy that is characteristic of the
spectral collocation methods.

\FloatBarrier

%


\section*{Results}
\label{sec3:results}

In this section, we analyze the \mymethod{} performance from various
perspectives.
First, we assess its accuracy and convergence by investigating
its behavior when approximating a low-rank, high-frequency analytical
solution.
The high-frequency test examines \mymethod{}'s performance on a grid
density with $300^7\approx2.2\, 10^{17}$ nodes that can only be explored by low-rank
solvers since a full format traditional solver, even at minimal
complexity would require millions of years on modern supercomputers.
Next, we present a six-dimensional test case that triggers the
\advection{} instabilities when we approximate a high-gradient
solution that challenges the low-rank approximation on moderate to
high-density grids. It tests the robustness of \mymethod{} in handling such instabilities.
Finally, we complete our investigation with two low-dimensional test cases
that challenge the non-oscillatory performance of
the \mymethod{} solver: the ``traveling bump''
benchmark~\cite{Funaro:1997,blue-book} and the Hughes double-layer
benchmark~\cite{Hughes-Mallet:1986,Brooks-Hughes:1982}.

\subsection*{Global Accuracy and Computational Cost}
\label{subsec:accuracy}

Testing numerical algorithms with manufactured analytical solutions
is a cornerstone practice in numerical analysis.
It provides a reliable way to create benchmark solutions for validating numerical solver\cite{manufactored-masa}.
In this case, we construct an ``exact'' solution by choosing a suitable expression,
and we determine the corresponding right-hand side analytically.
To this end, we apply our \mymethod{} solver to
problem~\eqref{eq:pblm:time:A}-\eqref{eq:pblm:time:C} on
$\Omega\times[0,T]$ with $\Omega=[-1,1]^6$ and $T=1$.
For convenience, we rescale the time interval to match the interval
$[-1,1]$.
We set $\varepsilon=10^{-4}$,
\begin{align*}
  \bm{\beta}(x_1,x_2,\ldots,x_6) &= ( \sin(\pi x_{1}), \sin(\pi x_{2}), \ldots, \sin(\pi x_{6}) )^T,\\[0.5em]
  \rho      (x_1,x_2,\ldots,x_6) &= 
  \prod_{1\leq i\leq 6}\left[\frac{1}{2}-\left(x_i+\frac{11}{10}\right)\left(x_i-\frac{1}{2}\right)
    \left(x_i+\frac{1}{2}\right)\left(x_i-\frac{11}{10}\right)\right].
\end{align*}
The right-hand side $b$ and the Dirichlet boundary conditions are
determined from the proposed analytical solution:
\begin{align*}
  \fs(\xv,\ts) = e^{\ts}
  \cdot\prod_{1\leq\is \leq 6}\left( \sin(\pi\xss{i}) + 0.3\sin(80\pi\xss{i}) \right).
\end{align*}
Here, we test the space-time formulation.
In this setting, the initial condition is naturally incorporated as a 
boundary condition at $t = -1$, while at $t = 1$, the numerical solution 
is a result of the numerical scheme.
%
The analytical solution $\fs(\xv,\ts)$ incorporates two spatial
components: a low-frequency one, i.e. $\sin(\pi\cdot)$, that spans the
computational domain, and a high-frequency one,
i.e. $\sin(80\pi\cdot)$, that introduces features at approximately
$1/100$ of the domain size.
Such multi-scale phenomena are commonly encountered in practical
applications~\cite{Glowinski-Pan:2022,Wienke:2005}.

The results are presented in Figure~\ref{fig:transport-convergence}.
Each ``dot'' represents an approximation of the solution of
problem~\eqref{eq:pblm:time:A} for a given polynomial degree
(associated with the number of degrees of freedom).
The convergence tests are based on the following $L^2$-type error
norm:
\begin{align*}
  \mathrm{error} = \frac{\norm{\fs-\fsn}{\calF}}{\norm{\fs}{\calF}},
\end{align*}
where $\norm{,\cdot,}{\calF}$ denotes the Frobenius norm and $\fsn$
represents the polynomial approximation of degree $\ns$ of the solution
$\fs$.

\newcommand{\PANEL}[1]{$\mathbf{(#1)}$} Panel \PANEL{b} shows that we need a grid of approximately $256^7$ size to capture the two dominant
frequencies of the solution.
From that point onwards, we can appreciate the exponential  
decay of the
error characteristic of the spectral methods.
The solver achieves an $\mathcal{O}(10^{-11})$ approximation error for
a grid size having a total number of $300^7$ nodes (including time).
Panel \PANEL{a} shows that \mymethod{} achieves the desired accuracy
in under three minutes on an Intel i9 laptop computer with 64 GBytes
RAM.
We put this number into perspective theorizing a full format 
Super Spectral Collocation solver with
minimal quadratic complexity $((\#\mathsf{dofs})^7)^2$, labeled 
\textit{full format solver}, \cite{Funaro:2002,Funaro:1997}.
We estimate that 
a modern supercomputer, characterized by a peak performance of
$2.7^{18}$ floating point operations per second\cite{TOP500},
would require approximately fifty million years to run the full format 
solver on the $256^7$ grid.
The computational cost for the first grid that exhibits spectral
accuracy, the one with $300^7$ total degrees of freedom, would approximately be 300 million years.
These two timescales respectively compare with the geological age of
the Himalayas and the Dolomites.
These results are a consequence of the fact that \mymethod{}
can capture a low-rank solution of a maximum rank of 9, resulting in
a compression ratio of
$(7\cdot9^{2}\cdot300 ) /{300^7}=\mathcal{O}(10^{-12})$. 
The corresponding memory requirement for the 
full format solver exceeds one exabyte, indicated by the arrow 
in Figure~\ref{fig:transport-convergence}-\PANEL{a}. 
One exabyte compares at best, or exceeds, the storage capacity of 
todays largest supercomputers \cite{TOP500}.

\begin{figure}[!t]
  \centering
  \includegraphics[width=1.\textwidth]{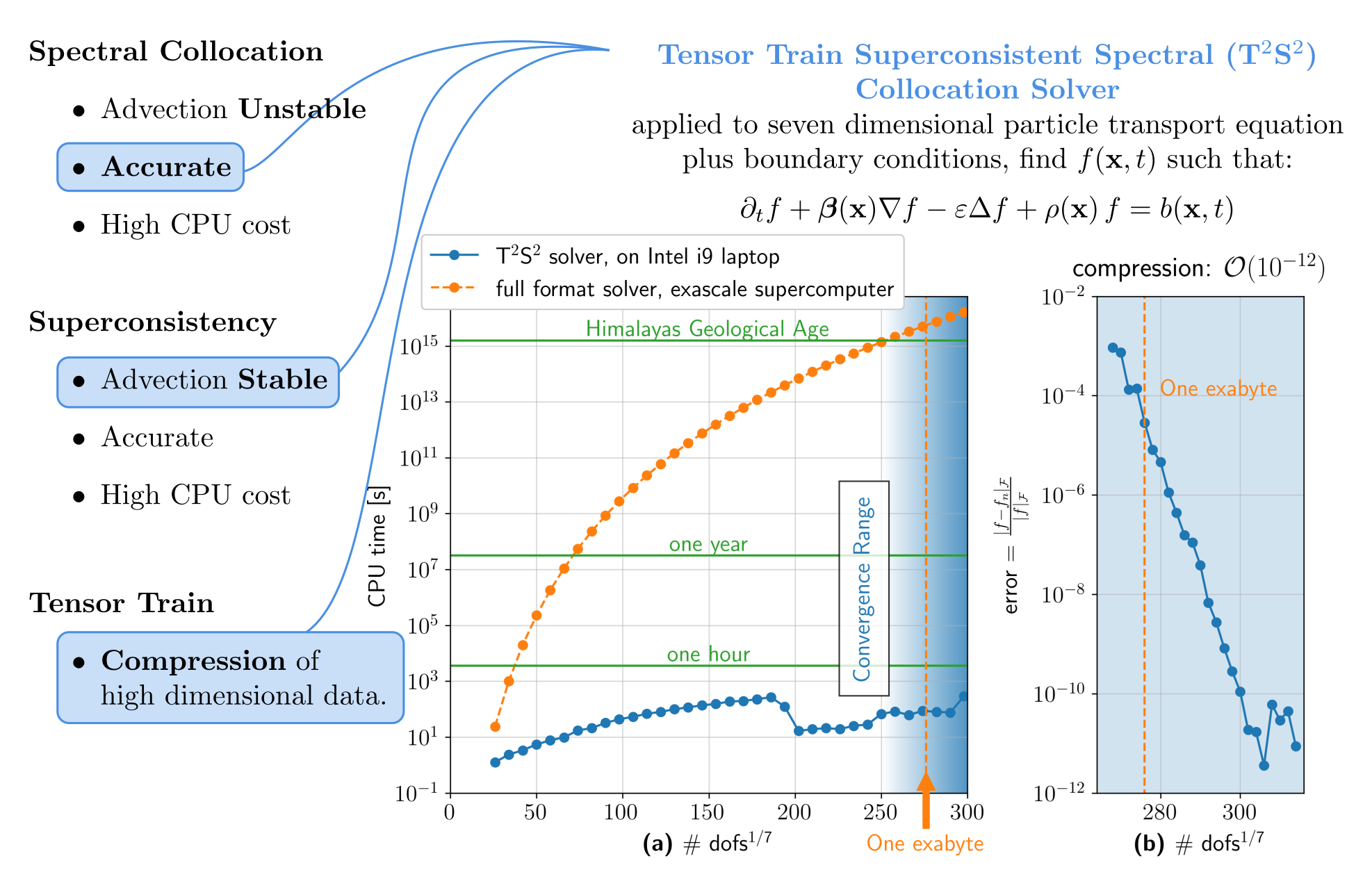}
  \caption{
    This figure showcases the key components of \mymethod{}  
    and its performance in solving the manufactured transport problem.
    Panel \PANEL{a} illustrates \mymethod{}'s compression
    capabilities, which enable the real-time solution of a
    seven-dimensional space-plus-time problem.
    To convey to the reader the scale of the performance, we provide 
    a comparison with the estimated computational time for a full format solver 
    with minimal numerical complexity~\cite{Funaro:1997,Canuto-Hussaini-Quarteroni-Zang:2007}, 
    whose cost scales as $\mathcal{O}((\#\mathsf{ndofs})^7)^2$. 
    Assuming such a full format solver runs on 
    a $2.7$ exaFLOPS nominal supercomputer, the projected 
    execution time would be on the order of fifty million years, comparable to the geological 
    age of the Himalayas. 
    An arrow is used to indicate the ``One exabyte'' memory for the full format solver.
    Panel \PANEL{b} shows the convergence behavior.
    The error decreases exponentially, eventually reaching machine
    precision.}
  \label{fig:transport-convergence}
\end{figure}

\FloatBarrier

\begin{figure}[t!]
  \centering
  \includegraphics[width=.5\textwidth]{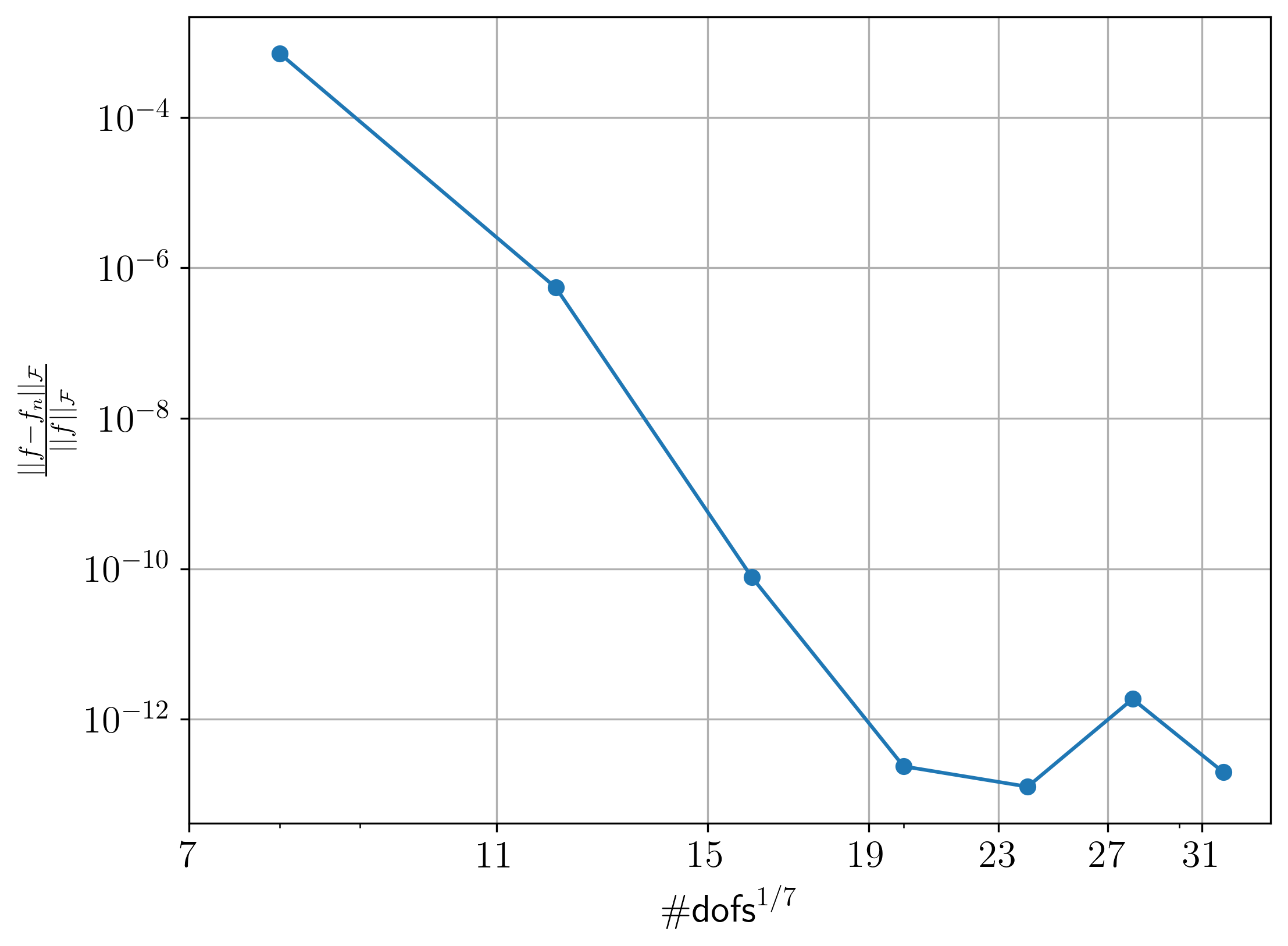}
  \caption{This convergence curve shows the exponential decay of the
    approximation error of the numerical solution versus the number of
    degrees of freedom in each direction from the application of
    \mymethod{} to a seven-dimensional convection-diffusion problem
    with velocity field $\bm{\beta}=(1,1,\ldots,1)^T$ and diffusion
    coefficient $\varepsilon=10^{-6}$.}
  \label{fig:transport-convergence-epsilon}
\end{figure}

\medskip
Figure~\ref{fig:transport-convergence-epsilon} shows the convergence analysis of \mymethod{} on a
six-dimensional, constant coefficient test case.
We set
$\varepsilon=10^{-6}$,
$\bm{\beta}$ as a six-dimensional unit vector and $\rho=0$.
The analytical solution is prescribed as
\begin{align*}
  \fs(\xv,t) = e^{-t} \cdot \prod_{ 1 \le i  \le 6}\left(\sin(\pi x_i) \right),
\end{align*}
from this, we derive the corresponding right-hand side and boundary
condition for the discrete model.
Figure~\ref{fig:transport-convergence-epsilon} presents the error
curve as a function of the number of degrees of freedom along one
dimension.
The results demonstrate that the method preserves spectral accuracy.

\FloatBarrier

\subsection*{Superconsistent Stabilization of Convection-Dominated Problems}
\label{subsection:stability}
Stabilization techniques for convection-diffusion problems have been
extensively studied in the
literature\cite{Hughes-Mallet:1986,Brooks-Hughes:1982}, mainly
focusing on low dimensional stationary convection-diffusion equation
with Dirichlet boundary conditions.
We consider a six dimensional equivalent:
\begin{align}
  \begin{split}
    -\varepsilon\Delta\fs + \bm{\beta}\cdot\nabla\fs &= \bs \phantom{\gs}\quad\textrm{in~}\Omega,\\[0.5em]
    \fs                                              &= 0 \phantom{\bs}\quad\textrm{on~}\partial\Omega.
  \end{split}
\label{eq:advection-diffusion}
\end{align}
This equation is a specialized case of our target problem
\eqref{eq:pblm:time:A}, which we analyze on the six-dimensional
hypercube domain $\Omega=[-1,1]^6$.
Conventional spectral approximations perform effectively in
diffusion-dominated regimes where $\varepsilon\Delta\fs$ predominates
over $\bm{\beta}\cdot\nabla\sf$. However, they exhibit spurious oscillations in
convection-dominated scenarios which are common in many practical
applications\cite{Glowinski-Pan:2022,Smith-Jensen:1989}.

\begin{figure}[!ht]
  \centering
  \includegraphics[width=0.8\textwidth]{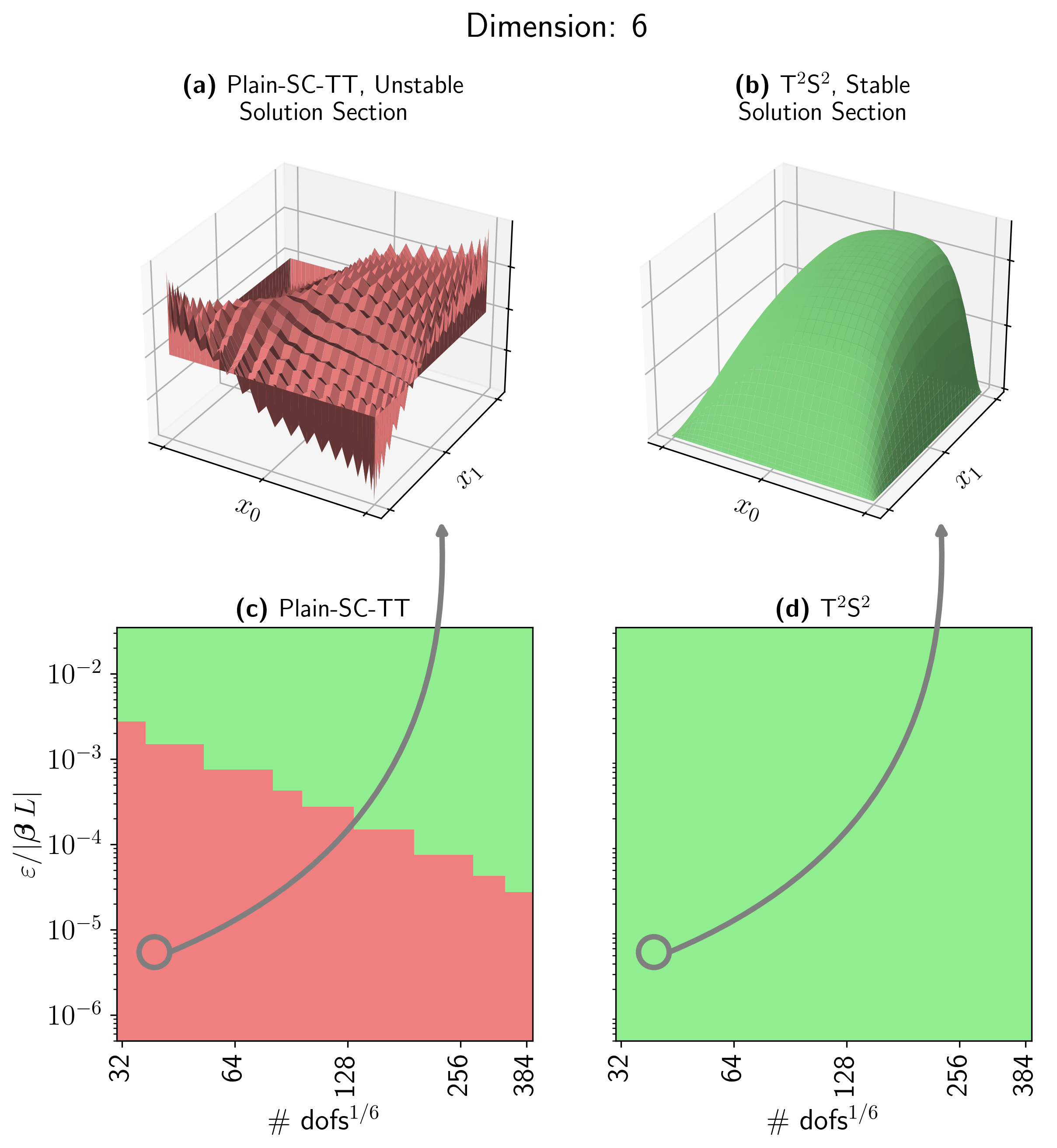}
  \caption{With $\bm\beta$ fixed as a six-dimensional unit vector, we
    conduct a detailed analysis of the solution's oscillatory behavior.
    Panel \PANEL{a} shows a two-dimensional cross-section of the
    solution to the six-dimensional problem
    (\ref{eq:advection-diffusion}) with $\varepsilon = 10^{-5}$ and
    polynomial degree 39, using collocation points at the zeros of $
    P'_n$. The solution presents 
    macroscopic spurious oscillations.
    The corresponding solution using superconsistent collocation
    points exhibits the expected regularity and characteristic
    boundary layers.
    Panels \PANEL{c} and \PANEL{d} present colormaps  
    where each pixel represents a discrete solution of
    problem \eqref{eq:advection-diffusion} for specific combinations
    of polynomial degree and \AdvRatio{}.
    Panel \PANEL{c} highlights regions of oscillatory solutions in
    red, demonstrating the limitations of \PLAINSC{}-TT.
    In contrast, panel \PANEL{d} illustrates how our method maintains
    non-oscillatory behavior across all values of \AdvRatio{} and
    polynomial degrees.
    }
  \label{fig:stabilitydomain}
\end{figure}

For our numerical experiments, we take $\bs=1$ and define
$\bm\beta$ as a vector of ones. The right-hand side is set as the
identity function, and homogeneous boundary conditions are applied.
In Figure~\ref{fig:stabilitydomain}, we vary the degree of the space
discretizations and $\varepsilon$. We compare \mymethod{} 
with the \PLAINSC{} stored in tensor train format, 
we label this method \PLAINSC{}-TT.
The results are displayed on colormaps \PANEL{c} and \PANEL{d}, 
where red pixels indicate an oscillatory solution;
green pixels represent a non-oscillatory solution.
We count the number of sign changes in the first derivative of the
numerical solution. We compare this number with the expected sing
changes in the regular solution. If the first count exceeds the second, the 
solution is marked as oscillatory.
Oscillations are evaluated on the mid-axis along the first direction,
at the center most degree of freedom.

Panels \PANEL{a}-\PANEL{b} of Figure \ref{fig:stabilitydomain} show a
two-dimensional slice of the six-dimensional solution computed with
$\varepsilon=10^{-5}$ and 40 degrees of freedom.
Panel \PANEL{b} presents the non-oscillatory solution, which exhibits a
continuous profile with a well-resolved boundary layer at the domain
perimeter.
In contrast, panel \PANEL{a} shows the oscillatory solution, where
spurious oscillations dominate the numerical approximation, rendering
the spatial discretization ineffective.
Colormap in Figure \ref{fig:stabilitydomain}-\PANEL{c} shows that the
\PLAINSC{}-TT approach is unstable for a wide range of parameters.
The stable/unstable interface follows a precise law; diffusive terms
scale like $\mathcal{O}(\#\mathsf{dofs}^{4})$, while transport terms
scale like $\mathcal{O}(\#\mathsf{dofs}^{2})$, see
Section~``Methods''.

\medskip
It is, in principle, true that even for a very small value of
\AdvRatio{}, we can find a sufficiently refined mesh where the
diffusion term becomes numerically dominant, and the \PLAINSC{}-TT
approximation is stable.
This approach is unfortunately illusive, since it is always possible
to choose smaller values of \AdvRatio{} that result in spurious
oscillations in the solution.
Moreover, the computational capability of Tensor Train linear solvers
is remarkable but not infinite.
AMEn\cite{Dolgov-Savostyanov:2014}, an advanced TT linear system
solver requires solving local problems at each TT core iteration, with
each local direct solution having complexity $\mathcal{O}((r_{\ell-1},
n_{\ell},r_{\ell})^{3})$, where $r_{\ell-1}, r_{\ell}$ denote the TT
ranks and $n_{\ell}$ the mode size of the corresponding core.
This results in a challenging computational limit on the maximum
possible resolution determined by the hardware's RAM.
This is why \mymethod{} comes into play.
The super-consistent discretization is always unconditionally stable
regardless of the choice of the parameters; see
Figure~\ref{fig:stabilitydomain}-\PANEL{d}.

The boundary layer in Figure \ref{fig:stabilitydomain}-\PANEL{b},
characterized by a high gradient connecting the domain solution to the
zero boundary condition challenges the rank compression of the Tensor
Train format, as higher gradients typically result in a more dense
solution structure in terms of dominant patterns \cite{Kutz:2013}.

\FloatBarrier

\begin{figure}
  \centering
  \includegraphics[width=0.75\textwidth]{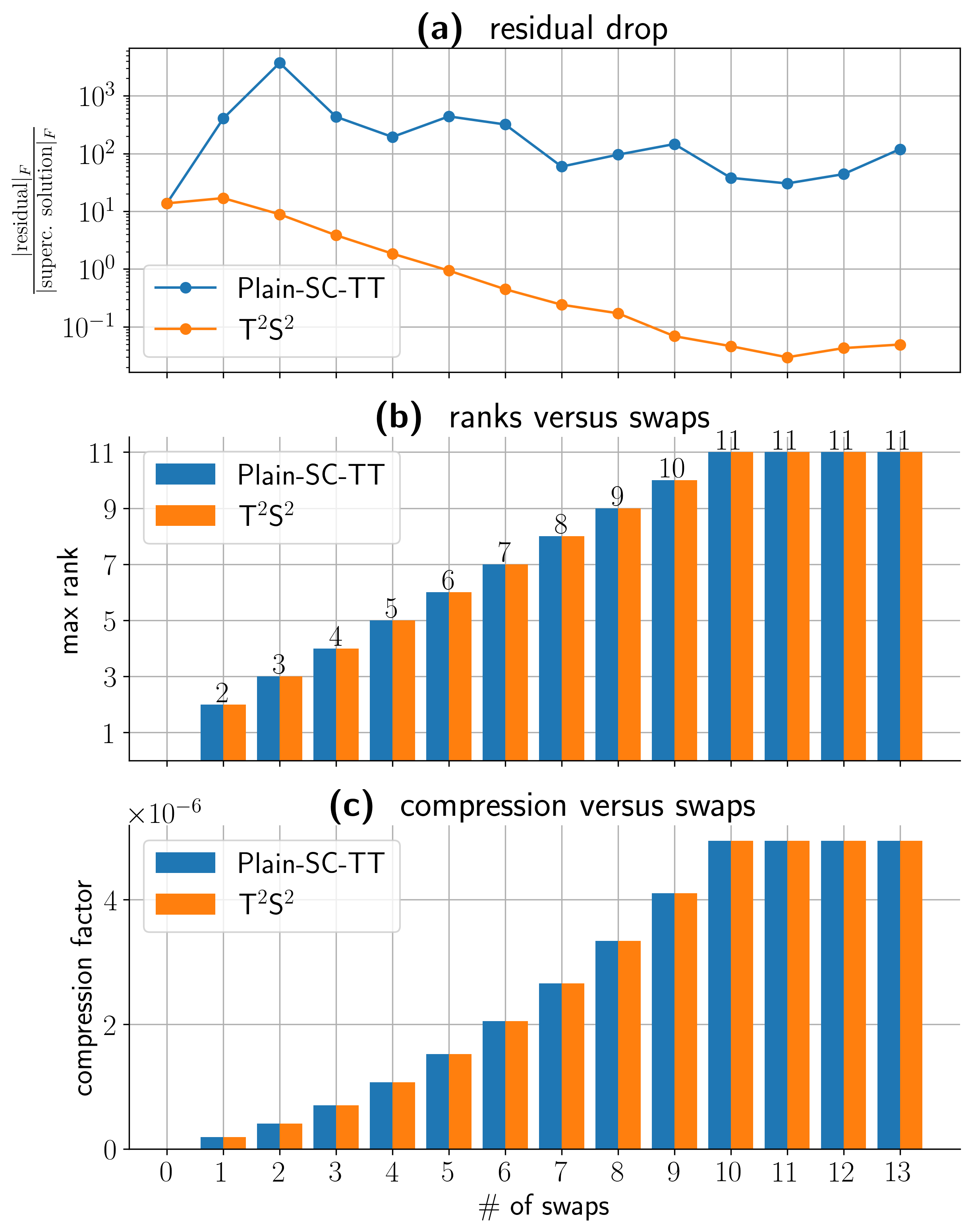}
  \caption{We present the residual decay, the ranks, and the
    compression factor against the number of AMEn swaps
    \cite{Dolgov-Savostyanov:2014} for our \mymethod{} solver and the
    \PLAINSC{}-TT collocation methods in the case of 6 dimensional
    problem \ref{eq:advection-diffusion}, degree 39,
    $\varepsilon=10^{-5}$ and $\bm{\beta}$ six-dimensional vector of
    ones.
    This problem is of special interest because it presents a boundary
    layer.
    The boundary layer is a high-gradient region of the solution that
    connects a boundary condition with the internal values.
    The approximation of such high gradients is generally challenging
    in the tensor train format.
    Panel \PANEL{a} reports how the residual drops against the number
    of swaps.
    The residual does not converge in the \PLAINSC{}-TT case, while in
    \mymethod{} case, the residual drops with a regular gradient.
    The final value of the residual $\mathcal{O}(10^{-2})$ is a
    consequence of the overall bad conditioning of the problem.
    Panel \PANEL{b} shows that the solver increases the rank of the
    solution at each swap until it finds a stationary point at swap 11
    for a corresponding rank equal to 11.
    After the 11th swap, the solver does not find any further
    information.
    The compression factor is in the range of $\mathcal{O}(10^{-6})$, \PANEL{c}.
    The computational cost on an i9 Intel laptop with 64 Gigabytes of RAM 
    is in the range of seconds.}
  \label{fig:residual-drop}
\end{figure}

In Figure \ref{fig:residual-drop}, we relate the 
linear solver ``iterations''
to the Tensor Train ranks, compression factors, 
and residuals for \mymethod{} and Plain-SC-TT. 
The Tensor train solver is AMEn\cite{Dolgov-Savostyanov:2014},
ant the ``iterations'' are addressed as swaps in the tensor train nomenclature.s
Panel \PANEL{b} illustrates that the tensor-train solver adjusts the
solution rank with each swap, reaching a stationary point at swap $11$
with a rank of $11$.
At this point, both our method and the \PLAINSC{}-TT approach reach their
resolution limit, and beyond it, they are unable to extract further
information.
Maximum rank $11$ corresponds to a compression factor around
$\mathcal{O}(10^{-6})$, cf. the details in panels \PANEL{b} and \PANEL{c}.
Although \mymethod{} and \PLAINSC{}-TT achieve a similar compression
ratios, their convergence properties diverge substantially.
The residual analysis in panel \PANEL{c} reveals that \mymethod{}
maintains consistent convergence with regular residual decay while
\PLAINSC{}-TT exhibits non-convergent behavior characterized by irregular
residual patterns.
The final residual value of $\mathcal{O}(10^{-2})$ reflects the
problem's poor conditioning
\cite{Funaro:2002,Canuto-Hussaini-Quarteroni-Zang:2007}.
Moreover, the computational cost of \mymethod{} for solving a problem
with $\#\mathsf{dofs}=\mathcal{O}(10^{9})$ on an Intel i9 laptop
computer is in the range of seconds.

\FloatBarrier

\subsection*{Low-Dimensional Numerical Challenges}
\label{subsection:low-dimensional-challenges}

\begin{figure}[!ht]
  \centering
  \includegraphics[width=0.85\textwidth]{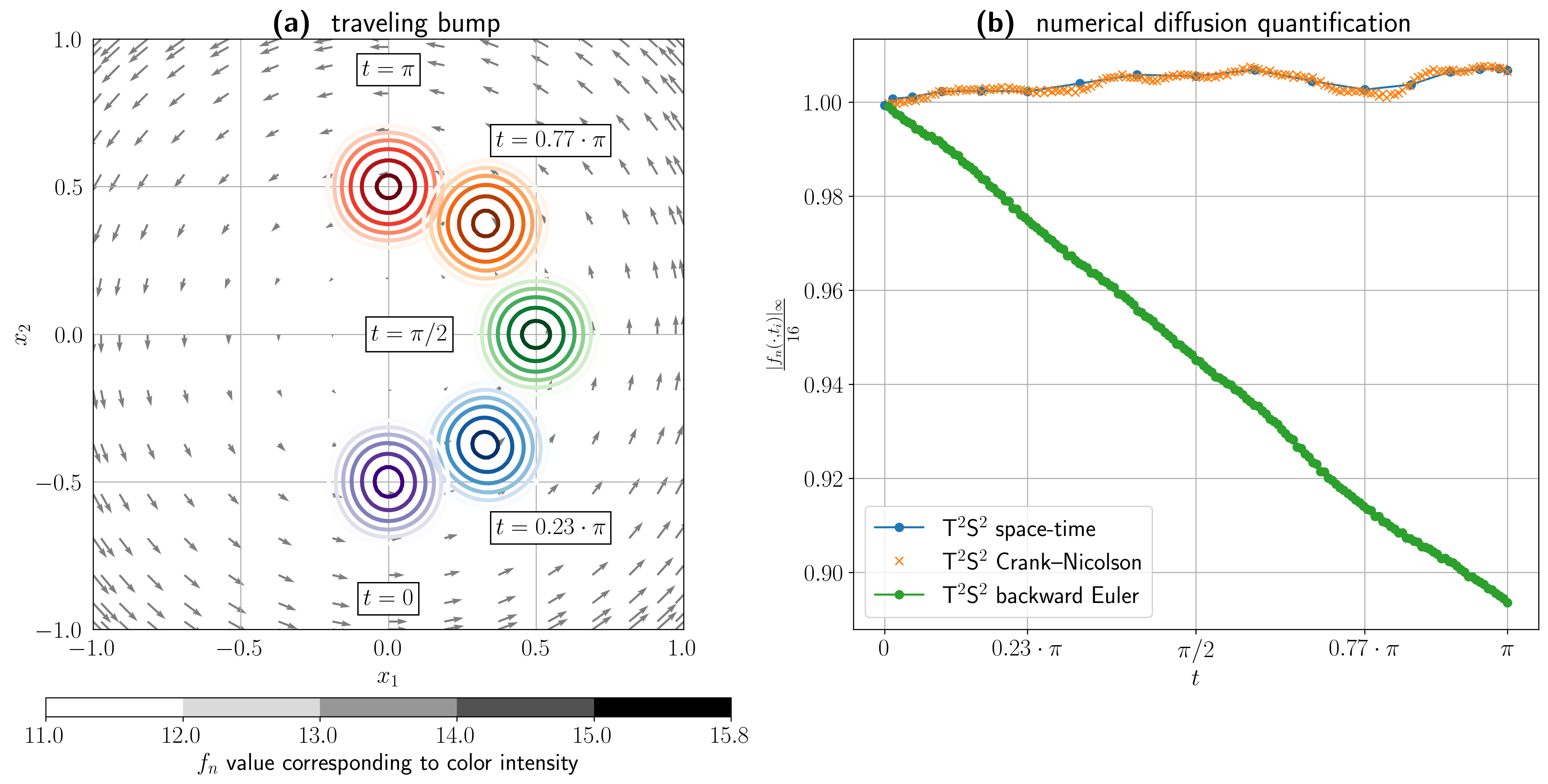}
  \caption{Results for the traveling ``bump'' test case.
    Panel \PANEL{a} shows contour lines for the ``bump'' at different
    time steps.
    Different colors distinguish the different instants; the shades
    distinguish the different function values.
    The different levels are reported in the gray color bar.
    Panel \PANEL{b} shows the time evolution of the solution peak
    values for three \mymethod{} method's implementations using the
    implicit backward Euler scheme, the semi-implicit Crank-Nicolson
    scheme, and the space-time formulation. A video of the traveling bump is provided as Supplementary Materials.}
  \label{fig:bulb-motion}
\end{figure}

This section presents two low-dimensional test cases to 
assess the characteristic capabilities of \mymethod{}.
The first test is the traveling bump
benchmark~\cite{blue-book,Funaro:1997}.
The bump benchmark is a two-dimensional time-dependent problem where
the initial condition is known, but there is no a priori knowledge of
the evolution of the solution.
The test is based on a rotational velocity flow, that challenges the
Cartesian structure of the \mymethod{} grid.
Although we choose collocation points to be the zeros of $P_n$,
considering positive convection coefficients in both directions,
\mymethod{} nonetheless proves to be robust.
We solve problem~\eqref{eq:pblm:time:A}-\eqref{eq:pblm:time:C} with
$\bm{\beta}=(-x_2,x_1)^T$, $\varepsilon=10^{-6}$ and $\rho = 0$.
The space-time domain is given by
$\Omega\times\Ts=\left([-1,1]\times[-1,1]\right)\times[0,\pi]$.
The initial condition is:
\begin{align*}
  f(x_1,x_2,0) = 
  \left\{
  \begin{array}{ll}
    16\,\left(1-4\,\left(x_2-\frac{1}{2}\right)^2\right)^2\,\left(1-4\,x_1^2\right)^2 & 
    x_2 > 0,\ \frac{1}{2} \leq x_1 <  \frac{1}{2},
    \\
    0 & \mathrm{elsewhere.}
  \end{array}
  \right.
\end{align*}
The polynomial degree of the basis functions is $16.$
Figure~\ref{fig:bulb-motion} presents the results of this calculation,
in the case of the spacetime formulation.
Panel \PANEL{a} illustrates the evolution of the bump's position
within the domain, which describes a circular trajectory through
contour lines at various time snapshots, together with the prescribed
velocity field.
Theoretically, the configuration at time $\pi$ should be mirror
symmetric versus the initial condition.
Panel \PANEL{b} provides quantitative accuracy analysis through the
normalized maximum norm $\lVert\fss{n}(\cdot,t_i)\rVert_\infty/16$,
computed on a $128\times128$ uniform grid.
The results demonstrate that \mymethod{} preserves the solution
structure with minimal oscillation behavior.
In this panel we superimpose the result for the spacetime formulation
and the time marching schemes: backward Euler and Crank-Nicolson.
The implicit backward Euler scheme reads:
\begin{align}
  f^{k+1} - \Delta\ts\Big( \varepsilon \Delta f^{k+1} - \bm{\beta}\cdot \nabla f^{k+1} \Big) =
  f^{k},
  \label{eq:ttss-backward-Euler}
\end{align} 
and the semi-implicit Crank-Nicolson scheme is:
\begin{align}
  f^{k+1} - \frac{\Delta\ts}{2}\Big( \varepsilon \Delta f^{k+1} - \bm{\beta}\cdot \nabla f^{k+1} \Big) =
  f^{k} + \frac{\Delta\ts}{2} ( \varepsilon \Delta f^{k} - \bm{\beta}\cdot \nabla f^{k} ).
  \label{eq:ttss-Crank-Nicolson}
\end{align} 
We choose a time step $\Delta t = \pi/160$ in both cases.
Since the implicit backward Euler scheme is only first-order accurate,
the solution approximation shows significative numerical diffusion as
the peak's value is reduced by roughly 10\%.
The Crank-Nicolson scheme is particularly well-suited for this problem
due to its reduced numerical dissipation compared to fully implicit
schemes.
This property is critical as it better preserves the system's physical
quantities, which is a key characteristic in many applicative problems.

\FloatBarrier

\begin{figure}[!ht]
  \centering
  \includegraphics[width=.5\textwidth]{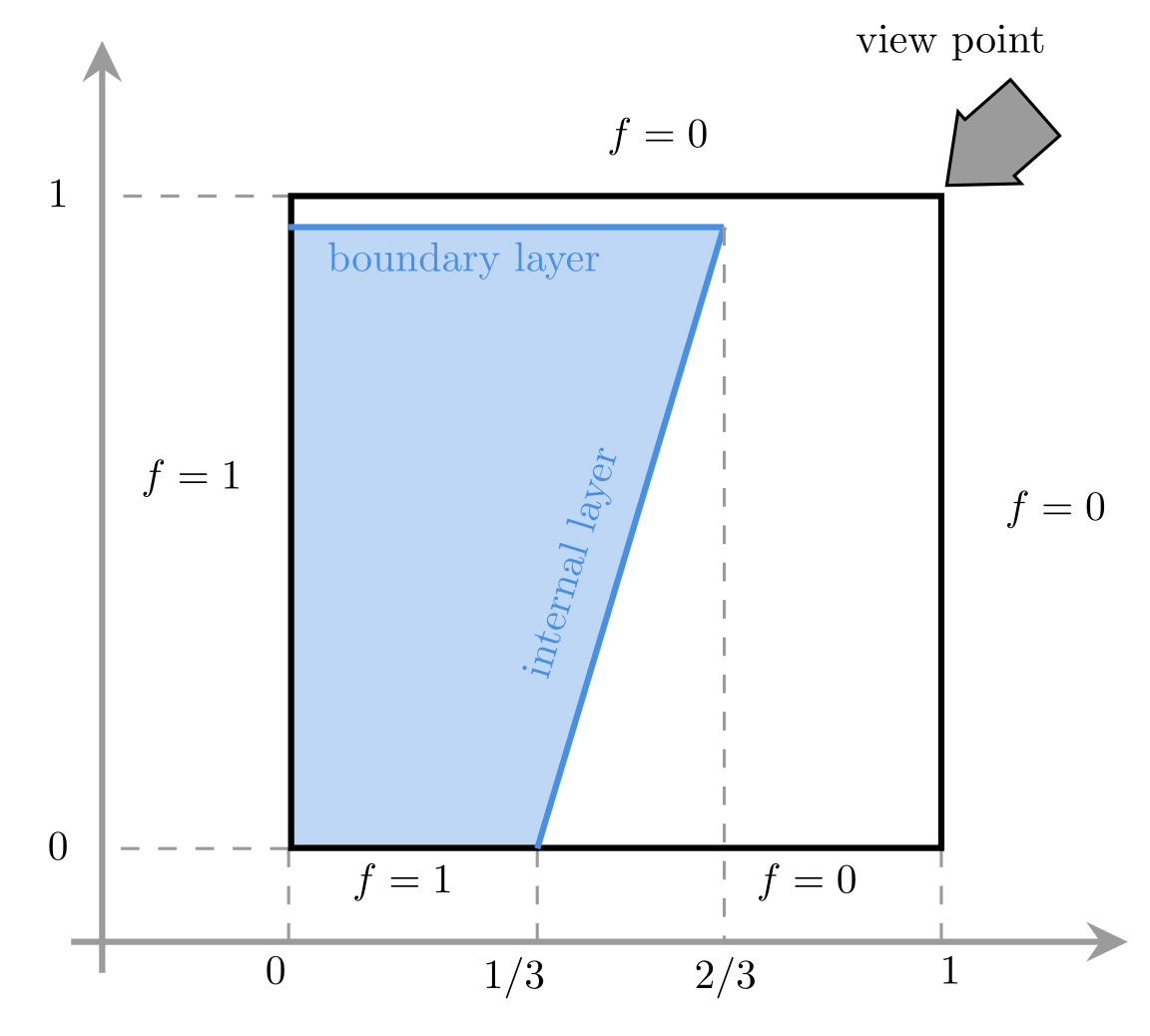}
  \caption{This picture gives a schematic overview of the 
  convection-diffusion test case as proposed by Hughes and
  coauthors\cite{Hughes-Mallet:1986,Brooks-Hughes:1982}, 
  highlighting the boundary and internal layers of the solution.}
  \label{fig:hughes-test-sketch}
\end{figure}

\begin{figure}[!t]
  \begin{center}
    \includegraphics[width=0.8\textwidth]{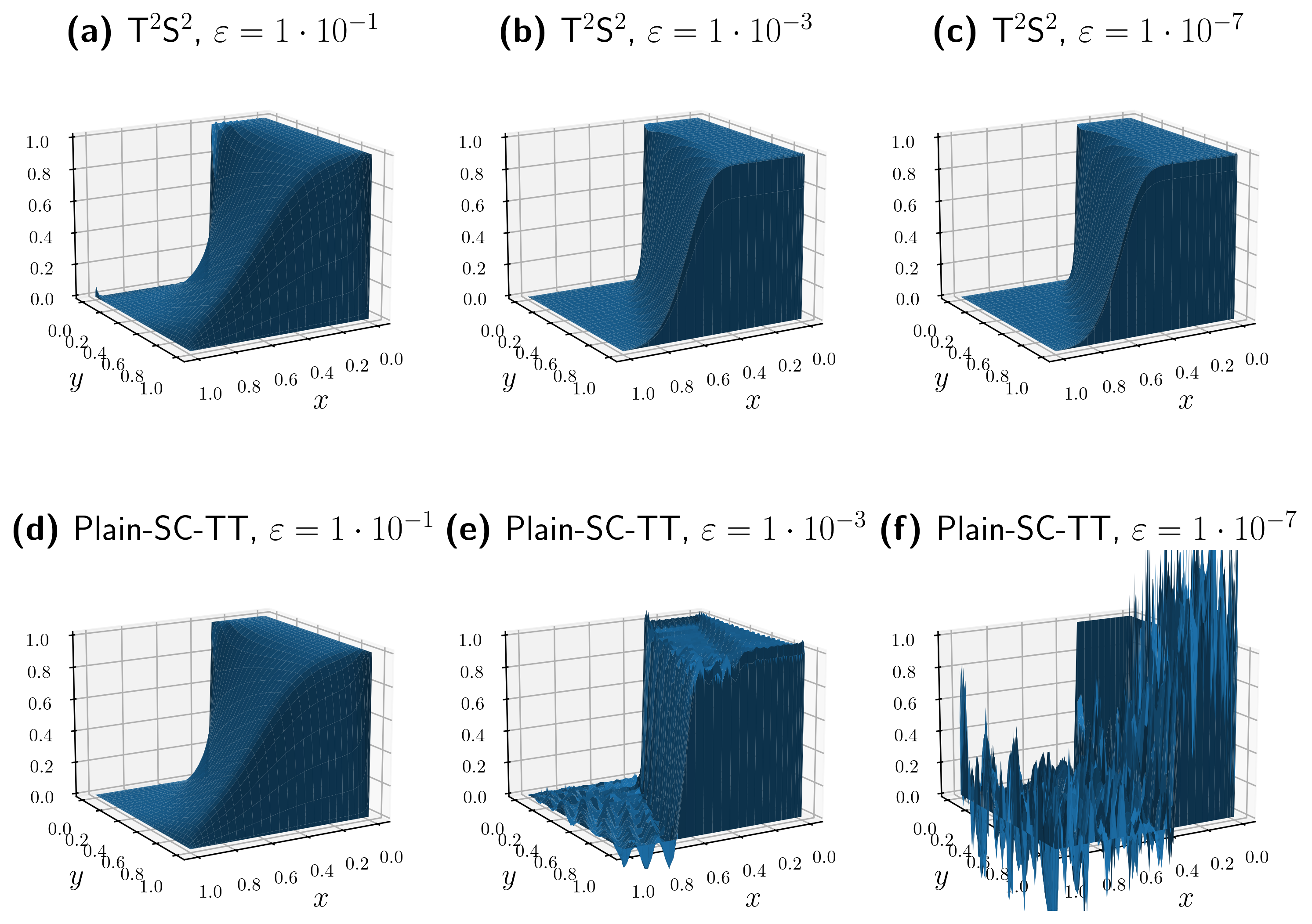}
  \end{center}
  \caption{This picture represents the solutions of the test proposed
    by Hughes and coauthors
    \cite{Hughes-Mallet:1986,Brooks-Hughes:1982} for various values
    of the parameter $\varepsilon$ and for the \PLAINSC{}-TT collocation
    method versus \mymethod{}.
    The solution degree is 63. While \PLAINSC{}-TT shows spurious oscillations at $\varepsilon=10^{-3}$, 
    \mymethod{} is stable regardless of the value of $\varepsilon$.
  \label{fig:hughes-test}}
\end{figure}

The second test case, originally proposed by Tom Hughes and
coauthors~\cite{Hughes-Mallet:1986,Brooks-Hughes:1982}, evaluates the method's ability to resolve both boundary and internal layers. We solve
problem \eqref{eq:advection-diffusion} in two dimensions with $\bs=0$
and decreasing values of $\varepsilon$.
The discontinuous boundary at $x_2=0$, which is transported in the
interior by the convection field $\bm{\beta}=(1,3)^T$, generates the
internal layer on the bottom side of the computational domain.
The boundary layer is determined by setting a homogeneous boundary
condition on the top side of the computational domain.
The test case specifics are sketched in
Figure~\ref{fig:hughes-test-sketch}.
The boundary layer connects the solution with the boundary conditions
and the internal layer connects two flat regions of the solution,
i.e., $f=1$ (on the left), and $f=0$ (on the right).
The reader will notice that the smaller is $\varepsilon$, 
the higher is the gradient of the internal layer; see
Figures~\ref{fig:hughes-test}-\PANEL{a}, \ref{fig:hughes-test}-\PANEL{b},
\ref{fig:hughes-test}-\PANEL{c}.
 \mymethod{} method confirms its oscillatory-free behavior. In contrast, \PLAINSC{} method exhibits spurious numerical oscillations, Figures~\ref{fig:hughes-test}-\PANEL{d}, \ref{fig:hughes-test}-\PANEL{e},
\ref{fig:hughes-test}-\PANEL{f}, which
start being visible for a relatively big value of
$\varepsilon=10^{-3}$, and fails to resolve the problem accurately.



\clearpage
\section*{Discussion}
\label{sec4:discussion}

\mymethod{} effectively  integrates the strengths of 
its core components: Spectral Collocation for accuracy, 
Superconsistency for stability and Tensor Train format as a low-rank 
linear algebra framework. 

In Section ``Results'' \ref{subsec:accuracy}, we demostrated that \mymethod{} reduces 
the approximation error to the range of the machine precision, 
with the error decreasing sharply when plotted on a logarithmic scale. 
Despite using a low-rank solution, the grid density is such
that traditional solvers---even at theoretically optimal efficiency--- 
would require millions of years on the current world’s most powerful architectures.\mymethod{} preserves spectral accuracy in 
high-dimensional, high dense grid scenarios.

In the section about stability, we experimentally tested 
the solver's stability using a six-dimensional stationary test 
case characterized by a boundary layer solution, i.e., a 
solution characterized by a high gradient region connecting the 
domain interior to the homogeneous boundary conditions. 
This test not only confirms the unconditional 
stability of the method versus the parameters choice, 
but also shows that \mymethod{} can approximate high 
dimensional/high gradient solutions with a 
rank in the range of 10.


\mymethod{} achieves this without requiring any intermediate approximation layers, 
either at the continuous or discrete level: the discrete version of problem~\eqref{eq:pblm:time:A} 
fits naturally into the Tensor Train structure. Moreover, \mymethod{} 
differs fundamentally from low-rank data-driven approaches, such as machine 
learning and reduced-order modeling\cite{Kutz-Brunton-Brunton-Proctor:2016,Hesthaven-Rozza-Stamm:2015}, 
as it does not require any high-fidelity offline computations or training phases. 
As discussed in ``Results'' \ref{subsec:accuracy}, we remarked that a high-fidelity 
computation in seven dimensions results in prohibitive computational
costs, even for a relatively coarse grid.

Following the characterization of low-rank methods proposed in
the comprehensive review by Einkemmer et al.\cite{Einkemmer:2024}, 
a time discretization of \mymethod{}, such as the backward Euler and Crank-Nicholson schemes, 
falls into the category of the \emph{Step And Truncate} (SAT) methods. 
SAT methods offer the flexibility to incorporate \mymethod{} into existing
time marching schemes, while spacetime \mymethod{} performs 
the compression also along the time dimension. 

Problem (\ref{eq:pblm:time:A}) serves as the basis for a wide
range of solvers in fields such as radiation transport physics, kinetic
theory, and computational fluid dynamics.
As such, \mymethod{} is a tool for a wide range of applications in modeling high-dimensional
transport phenomena.
Technological advances are driving the development of next-generation
nuclear technologies, including small modular reactors, advanced
modular reactors, microreactors, and nuclear fusion systems.
These innovations demand large-scale neutron transport and fluid dynamics
calculations for reactor design and throughout the nuclear fuel cycle,
including nuclear safeguards, nuclear non-proliferation, and security needs
\cite{ecp2025reports,iaea2025aris,international2024iaea,Geist-Santi-Swinhoe:2024,Cavallini-etal:2023}.
We expect \mymethod{} to be a key numerical solver in addressing
these complex computational challenges.

Beyond nuclear applications, \mymethod{} has the potential to impact a
wide range of fields requiring efficient, high-dimensional transport
computations.
Its ability to solve such problems on standard computing hardware
makes it particularly valuable for practical applications with limited
computational resources.
Future directions include extending the \mymethod{}
architecture to run efficiently on supercomputers, expanding
its applicability to large-scale scientific and engineering problems.


\section*{Methods}
\label{sec5:methods}

In this section, we describe the Tensor Train format
a fundamental 
building block for \mymethod{}, alongside 
Spectral Collocation Methods \cite{Funaro:1997}, and Superconsistency \cite{Funaro:2002}.


\subsection*{Tensor Train representation}

We first review some basic notation and definitions about
tensors~\cite{Kolda-Bader:2009} and, in particular, the tensor-train
format~\cite{Oseledets-Tyrtyshnikov:2010,Oseledets:2011}.
Then, we review some important basic facts about the tensor-train
decomposition.


\medskip
\paragraph{Notation and generalities.}

A $\DIM$-dimensional tensor
$\tA=\big\{\tA(\iss{1},\iss{2},\ldots,\iss{\DIM})\big\},\,\iss{\ell}=1,2,\ldots,\nss{\ell},\,\ell=1,2,\ldots,\DIM$
is a multi-dimensional array whose entries are indexed by $\DIM$
indices.
We denote $\DIM$-dimensional tensors using calligraphic, upper
case fonts, such as $\tA,\tB,\tC$.
For scalars, vectors, and matrices, which are formally $\DIM$-way
tensors with, respectively, a number of dimensions $\DIM=0,1,2$, we
adopt a different notation for clarity.
Precisely, we denote scalar quantities using normal, lower-case fonts
as in $\as,\bs,\cs$; vector fields using bold, lower-case fonts as in
$\av,\bv,\ccv$; matrices using normal, upper case fonts as in
$\As,\Bs,\Cs$.
Furthermore, we denote the entries of a vector $\av$ as $\as(i)$ and
the entries of a matrix $\As$ as $\As(i,j)$.
We may occasionally denote the entries of three-dimensional tensor as
$\tA(i,j,k)$ instead of $\tA(\iss{1},\iss{2},\iss{3})$, but we
prefer using the notation
$\tA(\iss{1},\iss{2},\ldots,\iss{\DIM})$ when $\DIM>3$.
For completeness, we anticipate that when we consider a
three-dimensional tensor that is also a core of a tensor train
decomposition (see below for the definition), the first and last
indices are a Greek letter, e.g., $\alpha, \beta$ etc, and the
intermediate index is always a Roman letter, e.g., $i,j$ etc.

Following a MATLAB-style notation, we denote free indices using the
symbol ``$:$''.
This notation makes it possible to represent the two important
structural concepts of \emph{fiber} and \emph{slice} of a tensor.
A fiber is a vector obtained by fixing all but one index, and a slice
is a matrix obtained by fixing all but two indices.
For example, a 3D tensor $\tA(i,j,k)$ has three different types of
fibers, i.e., $\tA(:,j,k)$, $\tA(i,:,k)$, $\tA(i,j,:)$, and three
different types of slices $\tA(i,:,:)$, $\tA(:,j,:)$, $\tA(:,:,k)$.
Finally, always according with the MATLAB-style notation, a repetition
of the free index symbol ``$:$'' denotes an index contraction.
For example, let $\As$, $\Bs$, $\Cs$ be three matrices with compatible
sizes such that $\Cs=\As\Bs$.
Then, the notation ``$\Cs(i,j)=\As(i,:)\Bs(:,j)$'' is equivalent to
$\Cs(i,j)=\sum_{\alpha}\As(i,\alpha)\Bs(\alpha,j)$, where the
summation over the index $\alpha$ is taken, as usual, over all the
elements of the $i$-th rows of $\As$ and $j$-th columns of $\Bs$.

\medskip
\paragraph{Tensor-train format.}

The tensor-train format provides an efficient representation of
high-dimensional tensors by decomposing them into a product of much
smaller, three-dimensional tensors called the ``TT cores''.
Specifically, a $\DIM$-dimensional tensor
$\tA\in\REAL^{\nss{1}\times\cdots\times\nss{\DIM}}$ is said to be in
tensor-train format if there exist $\DIM$ cores $\tAss{\ell} \in
\REAL^{\rss{\ell-1}\times\nss{\ell}\times\rss{\ell}}$, for
$\ell=1,\ldots,\DIM$, with $\rss{0}=\rss{\DIM}=1$, such that each
entry of $\tA$ can be computed as:
\begin{align}
  \tA(\iss{1},\ldots,\iss{\DIM}) =
  \sum_{\alpha_1=1}^{\rss{1}}\cdots\sum_{\alpha_{\DIM-1}=1}^{\rss{\DIM-1}}
  \tAss{1}(1,\iss{1},\alpha_1)\tAss{2}(\alpha_1,\iss{2},\alpha_2)\cdots\tAss{\DIM}(\alpha_{\DIM-1},\iss{\DIM},1).
  \label{eq:TT-format:elementwise}
\end{align}

In some sense, Equation~\eqref{eq:TT-format:elementwise} is a
multi-dimensional generalization of the (reduced) Singular Value
Decomposition (SVD) of matrices.
Note, however, that we do not require any specific orthogonality
property in the above definition (more details on the SVD connection
will be given in the next example and later in this section).

\begin{figure}[!t]
 \centering
 \includegraphics[width=0.8\textwidth]{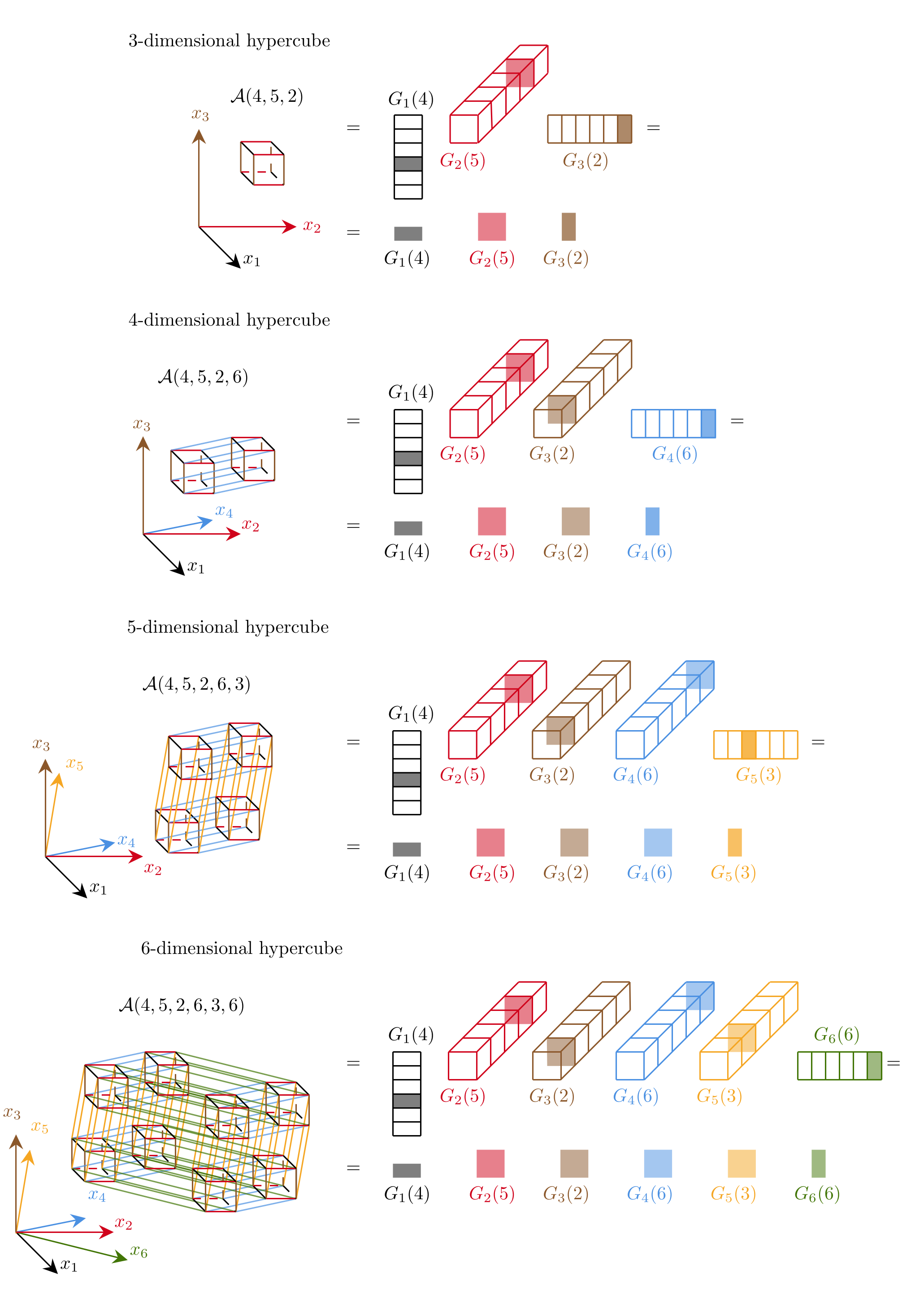}
 \caption{ In this picture we present the progression from a three to
   a six dimensional array in full and Tensor Train format.
 }
 \label{fig:tt-hypercube-six-dim}
\end{figure}

To fix ideas, consider the following two dimensional case.
A matrix $\As\in\REAL^{\nss{1}\times\nss{2}}$ is decomposed in the
product of two matrices $\Ass{1}\in\REAL^{\nss{1}\times\rss{1}}$ and
$\Ass{2}\in\REAL^{\rss{1}\times\nss{2}}$, so that
$\As(\iss{1},\iss{2})=\sum_{\alpha=1}^{\rss{1}}\Ass{1}(\iss{1},\alpha)\Ass{2}(\alpha,\iss{2})$.
We can easily construct such a decomposition through the SVD of $\As$,
which we write as usual as $\As=\Us\Sigma\Vs^T$, where
$\Us\in\REAL^{\nss{1}\times\rss{1}}$ and
$\Vs\in\REAL^{\nss{2}\times\rss{1}}$
are two orthogonal-by-column matrices, and
$\Sigma\in\REAL^{\rss{1}\times\rss{1}}$ is the diagonal matrix with
$\rss{1}$ strictly positive singular values, so that $\rss{1}$ is the
rank of matrix $\As$.
Then, take, for instance, $\Ass{1}=\Us$, $\Ass{2}=\Sigma\Vs^T$, or
$\Ass{1}=\Us\Sigma$, $\Ass{2}=\Vs^T$.
We can also use the tensor notation introduced above and specialized
to the case $\DIM=2$, so that we can write
$\tA(\iss{1},\iss{2})=\sum_{\alpha=1}^{\rss{1}}\tAss{1}(1,\iss{1},\alpha)\tAss{2}(\alpha,\iss{2},1)$,
where $\tAss{1}\in\REAL^{\rss{0}\times\nss{1}\times\rss{1}}$,
$\tAss{2}\in\REAL^{\rss{1}\times\nss{2}\times\rss{2}}$, with
$\rss{0}=\rss{2}=1$.
The decomposition process can be easily extended to $\DIM=3$ by
redefining $\tAss{2}$ as the intermediate core
$\tAss{2}\in\REAL^{\rss{1}\times\nss{2}\times\rss{2}}$ and introducing
a final core $\tAss{3}\in\REAL^{\rss{2}\times\nss{3}\times\rss{3}}$
with $\rss{3}=1$, so that
$\tA(\iss{1},\iss{2},\iss{3})
= \sum_{\alpha_1=1}^{\rss{1}}\sum_{\alpha_2=1}^{\rss{2}}
\tAss{1}(1,\iss{1},\alpha_1)\tAss{2}(\alpha_1,\iss{2},\alpha_2)\tAss{2}(\alpha_2,\iss{3},1)$.
This process can be continued by adding one new dimension at each step
as it is illustrated in Fig.~\ref{fig:tt-hypercube-six-dim} for the
construction of the TT representation of a tensor up to $\DIM=6$.

\paragraph{Remark.}~The computational effectiveness of the TT representation of
multi-dimensional arrays is related to controlling the core ranks.
To exemplify this fact, consider a six-dimensional array with mode
sizes $\ns=\nss{1}=\ldots=\nss{6}=10^3$ and assume that the ranks are
all equal to $\rs=\rss{1}=\ldots=\rss{d-1}=10$.
In such a case, the total storage required to represent such a tensor
in TT format is proportional to $\calO(\ns\rs^2d)=6\,10^5$ against a
full grid storage requirement of $\calO(\ns^d)=10^{18}$.

\medskip
The TT representation can equivalently be reformulated in a more
compact form by representing each core as a set of slice matrices
parametrized with the spatial index $\iss{\ell}$, e.g.,
$\Ass{\ell}(\iss{\ell})=\tAss{\ell}(:,\iss{\ell},:)\in\REAL^{\rss{\ell-1}\times\rss{\ell}}$,
yielding:
\begin{align}
  \tA(\iss{1},\ldots,\iss{\DIM}) = \Ass{1}(\iss{1})\Ass{2}(\iss{2})\cdots\Ass{\DIM}(\iss{\DIM}).
    \label{eq:TT-format:slice}
\end{align}
The internal sizes of the decomposition, i.e., the set of $(\DIM-1)$
integers $\{\rss{\ell}\}_{\ell=1}^{\DIM-1}$, are the the
\emph{TT-ranks}, or, simply, the \emph{ranks}, of the TT
decomposition.
A convenient upper bound on them is their maximum value, which we
denote as $\rs=\max_{\ell}\{\rss{\ell}\}$.
We also introduce an upper bound on the dimensions, denoted as
$\ns=\max_{\ell}\nss{\ell}$.
The TT-ranks determine the storage complexity of the representation,
which scales linearly with the dimension $\DIM$ as
$\mathcal{O}(\DIM\ns\rs^2)$, in contrast to the exponential growth
$\prod_{\ell=1}^{\DIM}\nss{\ell}=\calO(\ns^{\DIM})$ required for the
full tensor representation.
The linear dependence on $\DIM$ of the TT-representation of a tensor
is very convenient, reflecting the compression efficiency of such a
representation, whenever $\rs\ll\ns$.
In such a case, we say that \eqref{eq:TT-format:elementwise} and
\eqref{eq:TT-format:slice} are a \emph{low-rank representation} of
tensor $\tA$.

If we introduce the multi-index
$\is=(\iss{1},\iss{2},\ldots,\iss{\DIM})$, we can think that Equation
\eqref{eq:TT-format:elementwise} defines the TT format representation
for a sort of ``\emph{multi-dimensional vector}'', i.e., $\tA(\is)$.
It is also possible to adopt a similar \emph{TT format for linear
operators}, i.e., ``\emph{multi-dimensional matrices}'' $\calL(\is,\js)$ that depend on
two $\DIM$-dimensional multi-indices
$\is=(\iss{1},\iss{2},\ldots,\iss{\DIM})$ and
$\js=(\jss{1},\jss{2},\ldots,\jss{\DIM})$.
The formal definition reads as
\begin{align}
  &\calL\Big(\big(\iss{1},\iss{2},\ldots,\iss{d}\big),\,\big(\jss{1},\jss{2},\ldots,\jss{d}\big)\Big)
  =\calL\Big(\big(\iss{1},\jss{1}\big),\big(\iss{2},\jss{2}\big),\ldots,(\iss{d},\jss{d}\big)\Big)\nonumber\\[0.5em]
  &\hspace{2cm}
  =
  \sum_{\alpha_1=1}^{\rss{1}}\sum_{\alpha_2=1}^{\rss{2}}\ldots\sum_{\alpha_{\DIM-1}=1}^{\rss{\DIM-1}}
  \Lss{1}\big(1,(\iss{1},\jss{1}),\alpha_1\big)\,
  \Lss{2}\big(\alpha_2,(\iss{2},\jss{2}),\alpha_3\big)\, \ldots
  \Lss{d}\big(\alpha_{d-1},(\iss{d},\jss{d}),1\big),
  \label{eq:tt-matrix-format}
\end{align}
where the first equality follows from a convenient index permutation
to form the index pairs $(\iss{\ell},\jss{\ell})$,
$\ell=1,2,\ldots,d$.

The so-called ``\emph{core notation}'' gives a more intuitive insight
into the structure of the cores:
\begin{align*}
  \Ass{\ell}=
  \left[
    \begin{array}{cccc}
      \Ass{\ell}(1,:,1)            & \Ass{\ell}(1,:,2)            & \ldots & \Ass{\ell}(1,:,\rss{\ell})           \\
      \Ass{\ell}(2,:,1)            & \Ass{\ell}(2,:,2)            & \ldots & \Ass{\ell}(2,:,\rss{\ell})           \\
      \vdots                       & \vdots                       & \vdots & \vdots                               \\
      \Ass{\ell}(\rss{\ell-1},:,1) & \Ass{\ell}(\rss{\ell-1},:,2) & \ldots & \Ass{\ell}(\rss{\ell-1},:,\rss{\ell})\\
    \end{array}
    \right],\quad
\end{align*}
\begin{align*}
  \Lss{\ell}=
  \left[
    \begin{array}{cccc}
      \Lss{\ell}(1,:,:,1)            & \Lss{\ell}(1,:,:,2) & \ldots & \Lss{\ell}(1,:,:,\rss{\ell})           \\
      \Lss{\ell}(2,:,:,1)            & \Lss{\ell}(2,:,:,2) & \ldots & \Lss{\ell}(2,:,:,\rss{\ell})           \\
      \vdots                       & \vdots            & \vdots & \vdots                               \\
      \Lss{\ell}(\rss{\ell-1},:,:,1) & \Lss{\ell}(\rss{\ell-1},:,:,2) & \ldots & \Lss{\ell}(\rss{\ell-1},:,:,\rss{\ell})\\
    \end{array}
    \right].
\end{align*}
With such notation and using the definition of the strong Kroncker
product \cite{DeLauney-Seberry:1994} ``$\bowtie$'', we write the TT
decomposition in the more compact form:
\begin{align*}
  \mathcal{A} = \Ass{1}\bowtie\Ass{2}\bowtie\ldots\bowtie\Ass{d},\qquad
  \mathcal{L} = \Lss{1}\bowtie\Lss{2}\bowtie\ldots\bowtie\Lss{d}.
\end{align*}
We will make use of this mathematical
formalism to describe the different parts of \mymethod{}.

\medskip
\paragraph{Tensor-train (TT) decomposition}
The construction of an optimal tensor-train representation for a
$\DIM$-dimensional tensor can be achieved through the \TTSVD{}
algorithm~\cite{Oseledets:2011}.
This algorithm builds the TT decomposition through a systematic
application of singular value decompositions on specific matrix
representations of the tensor, known as \emph{unfoldings}.
For a $\DIM$-dimensional tensor $\tA$, we consider $\DIM-1$
``fundamental'' matrix unfoldings.
The $\ell$-th unfolding, denoted as $\As_\ell$, reorganizes the tensor
entries into a matrix of size
$(\prod_{k=1}^\ell\nss{k})\times(\prod_{k=\ell+1}^\DIM \nss{k})$
through the relation:
\begin{align}
  \As_\ell( \eta_\ell(\iss{1},\ldots,\iss{\ell}), \eta'_\ell(\iss{\ell+1},\ldots,\iss{\DIM})) 
  = \tA(\iss{1},\ldots,\iss{\DIM}),
\end{align}
where $\eta_\ell$ and $\eta'_\ell$ are two bijective mappings between
the respective index tuples, e.g., $(\iss{1},\ldots,\iss{\ell})$ and
$(\iss{\ell+1},\ldots,\iss{\DIM})$, often called ``long indices'', and
the row and column indices of matrix $\As_\ell$.
A very strong connection exists between the TT-ranks of tensor $\tA$
and the ranks of the unfolding matrices.
In fact, the rank of the $\ell$-th matrix unfolding is a lower bound
for the TT-rank $\rss{\ell}$, and a constructive proof shows that
\TTSVD{} can compute an \emph{exact} TT-representation with TT-ranks
$\rss{\ell}$ equal to the ranks of the matrix unfoldings
$\As_\ell$~\cite[Theorem~2.1]{Oseledets:2011}, i.e.,
$\rss{\ell}=\operatorname{rank}(\As_{\ell})$ for
$\ell=1,\ldots,\DIM-1$.

Of particular practical importance is the existence of low-rank
approximations.
The sequential nature of the \TTSVD{} algorithm ensures that
the TT ranks are optimal in the sense that they provide the best
possible approximation for a given accuracy threshold.
If we truncate the singular values in the \TTSVD{} algorithm below a
give threshold $\epsilon$, the resulting approximation error in the
Frobenius norm is bounded by $\epsilon\sqrt{\DIM-1}\|\tA\|_F$.
Such an approximate TT decomposition keeps the TT ranks minimal while
achieving the accuracy prescribed by $\epsilon$.
Therefore, when computational constraints or data compression
requirements need smaller ranks
$r_{\ell}'\leq\operatorname{rank}(\As_{\ell})$, a best approximation
$\tATT_{\text{best}}$ in the tensor-train format is guaranteed to
exist.

The \TTSVD{} algorithm efficiently provides a
quasi-optimal approximation $\tATT$ in the sense that:
\begin{align}
  \|\tA - \tATT\|_{\tF} \leq
  \sqrt{\DIM-1}
  \inf_{\tB^{\TT}\in\cTT(\rv'_\ell)} \|\tA - \tB^{\TT}\|_{\tF}.
\end{align}
Here, $\cTT(\rv'_\ell)$ denotes the manifold of tensor-train
representations with prescribed rank $\rv'_\ell =
(r'_1,\ldots,r'_{\DIM-1})$, and the multiplicative factor
$\sqrt{\DIM-1}$ represents the stability constant of the
algorithm~\cite{Oseledets-Tyrtyshnikov:2010}.
This quasi-optimality result establishes \TTSVD{} as a practical tool
for tensor approximation, offering a controlled trade-off between
computational efficiency and approximation accuracy.

The tensor-train format supports the efficient implementation of
fundamental multi-linear algebraic operations such as
scalar multiplication, tensor addition,
element-wise (Hadamard) products, and contractions with matrices and
vectors~\cite{Oseledets:2011}.
These operations are designed to preserve the tensor-train structure
but they generally lead to an increase in the TT ranks of the
resulting tensor.
A critical step is the so-called rounding procedure, denoted as
$\RNDG(\,\cdot\,)$, which ``recompresses'' the tensor by computing a
new TT representation with smaller ranks while controlling the
approximation error~\cite{Oseledets:2011}.
Given a tensor $\tA$ in TT-format, the ``rounded'' tensor $\RNDG(\tA)$
satisfies the inequality $\|\tA - \RNDG(\tA)\|_{\tF} \leq
\toll\|\tA\|_{\tF}$, where $\toll$ is a threshold parameter that makes
a trade-off between accuracy and computational costs possible.
On the one hand, smaller values of $\toll$ can preserve a better accuracy
but implies higher TT ranks, hence increasing the computational costs
and storage requirements.
On the other hand, larger values of $\toll$ lead to a better
compression with lower ranks, thus reducing computational complexity
at the expense of accuracy.
This balance between accuracy and efficiency must be carefully
calibrated according to the specific requirements of the application
at hand.

Finally, it is worth mentioning the cross-approximation
technique~\cite{Oseledets-Tyrtyshnikov:2010}, which allows us to
compute the TT representation of a tensor by sampling only a small subset of its entries.
This algorithm makes the calculation of a TT decomposition very
efficient since in many applications, mainly when dealing with
function-generated tensors, computing and storing the full tensor is
prohibitively expensive or impossible.
Cross-interpolation typically employs the \MAXVOL{}
algorithm\cite{Goreinov-Oseldets-etal:2010}, a sophisticated sampling
strategy that selectively chooses optimal rows and columns from the
tensor to construct its TT representation.
In our context, cross-interpolation proves particularly valuable for
efficiently representing the right-hand side term $\fs$, allowing us
to work with high-dimensional data while maintaining computational
tractability.




\printbibliography



\end{document}